\documentclass[a4paper,11pt]{article}
\usepackage[numbers]{natbib}
\usepackage{amsmath}
\numberwithin{equation}{section}
\usepackage{amssymb}
\usepackage{amsthm}
\usepackage{enumitem} 
\usepackage{graphicx}
\usepackage[tight]{subfigure}
\usepackage[font=small,labelfont=bf]{caption}
\usepackage{titlesec}
\usepackage[left=30mm,right=30mm,top=30mm,bottom=30mm]{geometry}

\usepackage{algorithm}
\usepackage[noend]{algpseudocode}
\usepackage{xpatch}
\makeatletter
\def\BState{\State\hskip-\ALG@thistlm}
\xpatchcmd{\algorithmic}{\itemsep\z@}{\itemsep=0.25ex}{}{}
\makeatother

\newtheorem{proposition}{Proposition}
\theoremstyle{definition}
\newtheorem{remark}{Remark}

\newcommand{\email}[1]{{#1}}
\renewcommand{\vec}[1]{\mathbf{#1}}
\newcommand{\x}{\vec{x}}
\newcommand{\f}{\vec{f}}
\newcommand{\Real}{\mathbb{R}}
\newcommand{\Sprocedure}{$\mathcal{S}$-procedure}
\newcommand{\expect}[2]{\langle #1 \rangle_{#2}}
\newcommand{\del}{\nabla}
 
\newcommand{\e}{\varepsilon}
\renewcommand{\phi}{\varphi}
\newcommand{\abs}[1]{\vert #1 \vert}
\newcommand{\norm}[1]{\| #1\|}

\newcommand{\bcolvec}[2]{\begin{bmatrix} #1 \\ #2 \end{bmatrix}}

\newcommand{\set}[1]{\lbrace #1\rbrace}

\title{Bounds for deterministic and stochastic dynamical systems using 
sum-of-squares optimization} 

\date{}

\author{
G. Fantuzzi\footnotemark[2] 
\and D. Goluskin\footnotemark[3]
\and D. Huang\footnotemark[4]
\and S.~I. Chernyshenko\footnotemark[2]
}

\begin{document}
%%%%%%%%%%%%%%%%%%%%%%%%%%%%%%%%%%%%%%%%%%%%%%%%%%%%%%%%%%%%%%%%%%%%%%%%%%%%%%%

%%%%%%%%%%%%%%%%%%%%%%%%%%%%%%%%%%%%%%%%%%%%%%%%%%%%%%%%%%%%%%%%%%%%%%%%%%%%%%%
% TITLE
%%%%%%%%%%%%%%%%%%%%%%%%%%%%%%%%%%%%%%%%%%%%%%%%%%%%%%%%%%%%%%%%%%%%%%%%%%%%%%%
\maketitle

\renewcommand{\thefootnote}{\fnsymbol{footnote}}
\footnotetext[2]{Department of Aeronautics, Imperial College London, South 
Kensington Campus, London, SW7 2AZ, United Kingdom 
(\email{gf910@ic.ac.uk}, 
\email{s.chernyshenko@imperial.ac.uk})}
\footnotetext[3]{Mathematics Department \& Center for the Study of Complex 
Systems, University of Michigan, Ann Arbor, MI, USA 
(\email{goluskin@umich.edu})}
\footnotetext[4]{School of Electrical Engineering, Southwest Jiaotong 
University, Chengdu, 610031, China (\email{elehd2012@gmail.com})}
\renewcommand{\thefootnote}{\arabic{footnote}}

%%%%%%%%%%%%%%%%%%%%%%%%%%%%%%%%%%%%%%%%%%%%%%%%%%%%%%%%%%%%%%%%%%%%%%%%%%%%%%%
% ABSTRACT
%%%%%%%%%%%%%%%%%%%%%%%%%%%%%%%%%%%%%%%%%%%%%%%%%%%%%%%%%%%%%%%%%%%%%%%%%%%%%%%
\begin{abstract}
We describe methods for proving upper and lower bounds on infinite-time averages 
in deterministic dynamical systems and on stationary expectations in stochastic 
systems. The dynamics and the quantities to be bounded are assumed to be 
polynomial functions of the state variables. The methods are computer-assisted, 
using sum-of-squares polynomials to formulate sufficient conditions that can be 
checked by semidefinite programming. In the deterministic case, we seek tight 
bounds that apply to particular local attractors. An obstacle to proving such 
bounds is that they do not hold globally; they are generally violated by 
trajectories starting outside the local basin of attraction. We describe two 
closely related ways past this obstacle: one that requires knowing a subset of 
the basin of attraction, and another that considers the zero-noise limit of the 
corresponding stochastic system. The bounding methods are illustrated using the 
van der Pol oscillator. We bound deterministic averages on the attracting limit 
cycle above and below to within 1\%, which requires a lower bound that does not 
hold for the unstable fixed point at the origin. We obtain similarly tight upper 
and lower bounds on stochastic expectations for a range of noise amplitudes. 
Limitations of our methods for certain types of deterministic systems are 
discussed, along with prospects for improvement.
\end{abstract}

%%%%%%%%%%%%%%%%%%%%%%%%%%%%%%%%%%%%%%%%%%%%%%%%%%%%%%%%%%%%%%%%%%%%%%%%%%%%%
%% MAIN TEXT
%%%%%%%%%%%%%%%%%%%%%%%%%%%%%%%%%%%%%%%%%%%%%%%%%%%%%%%%%%%%%%%%%%%%%%%%%%%%%

%%%%%%%%%%%%%%%%%%%%%%%% 
%    INTRODUCTION
%%%%%%%%%%%%%%%%%%%%%%%%
\section{Introduction}
\label{S:introduction}

In the study of dynamical systems with complicated and possibly chaotic 
dynamics, average quantities are often of more interest than any particular 
solution trajectory. This is partly because of the difficulty of computing a 
trajectory precisely and partly because average quantities are more important in 
many applications. For instance, one might seek time-averaged drag forces in a 
model of an oil pipeline or ensemble-averaged temperatures in a stochastic 
climate model. One way to estimate averages quantitatively is to integrate the 
dynamical system numerically and average over the resulting particular solution. 
Such direct numerical simulations are often straightforward, but the 
accuracy of the result is not guaranteed unless errors are rigorously 
controlled. Moreover, even a perfectly accurate solution does not give 
information about the different trajectories that result from different initial 
conditions. Finally, the computational cost of computing well-converged averages 
is often prohibitive, especially in systems that are high-dimensional or 
stochastic.

A complementary approach, which we pursue here, is to prove bounds on average 
quantities directly from the governing equations. An advantage over numerical 
integration is that bounds can be proven without knowing any solution 
trajectories. Furthermore, they can be proven for all possible trajectories at 
once, or for all trajectories within a given region of state space. On the other 
hand, it is generally difficult to prove bounds that are tight enough to give 
good estimates of the average quantities being bounded. 

The aim of this work is to develop methods for proving bounds that are tight, 
meaning that the upper and lower bounds are equal or nearly so. We assume that 
the quantity of interest can be described by a function $\phi(\x)$, where 
$\x(t)\in X$ is the state vector of a dynamical system, and we seek to bound 
averages of $\phi$. In deterministic systems, we consider averages over infinite 
time,
\begin{equation}
\label{E:TimeAvDef}
\overline{\phi(\x)} := \lim_{T\to\infty} \frac{1}{T} \int_{0}^{T} \phi[\x(t)] 
\,dt.
\end{equation}
If the above limit does not exist it can be replaced by limit superior or inferior for the upper and lower bound problems, respectively. The value of $\overline\phi$ depends in general on which trajectory $\x(t)$ is being averaged over. 
In stochastic systems with a stationary probability 
distribution $\rho(\x)$, we consider stationary ensemble averages,
\begin{equation}
\label{E:ExpectDef}
\expect{\phi(\x)}{} := \int_X \phi(\x) \rho(\x) \,d\x.
\end{equation}

One obstacle to proving tight bounds is, for reasons described shortly, the need 
to determine whether certain complicated expressions are sign-definite. As 
proposed in~\cite{Chernyshenko2014a}, this can be done systematically with 
computer assistance for finite-dimensional systems with polynomial dynamics. The 
main idea, described further in \S\ref{SS:SoSBoundsReview}, is to construct a 
polynomial whose non-negativity implies the desired bound. By the methods of 
sum-of-squares (SoS) programming~\cite{Parrilo2000,Parrilo2003}, 
a sufficient condition 
for this non-negativity can then be posed as a semidefinite program (SDP).

In deterministic systems there is a second obstacle to proving tight bounds. It 
is generally easiest to construct bounds that hold for all possible initial 
conditions. 
Sometimes this is desired, but other times one is interested in a particular 
local attractor, and bounds holding globally are not generally tight for 
averages over a local attractor. In this sense, the local bound is spoiled by 
any other invariant structure in the state space, such as another attractor or 
an unstable fixed point. 

We pursue two ways of obtaining tight bounds specific to a local attractor. The 
first is to enforce conditions that imply the bound only on an absorbing set 
around the attractor, thereby omitting other invariant structures. The second is 
to add noise to the system and prove bounds for ensemble averages in the 
vanishing noise limit. If the system is stochastically stable, 
then under certain conditions this limit will agree with the corresponding 
deterministic time average~\cite{Young2002}. 
Note that these ideas are applicable irrespectively of any special structure 
in the system (such as Hamiltonian), and can in principle
be applied to systems of arbitrarily high but finite dimension.

\begin{figure}[t!]
\centering
\subfigure[]{\label{F:vdpPhasePortrait}
\includegraphics[width=0.45\textwidth, trim=0cm 3.5cm 5cm 0cm]
{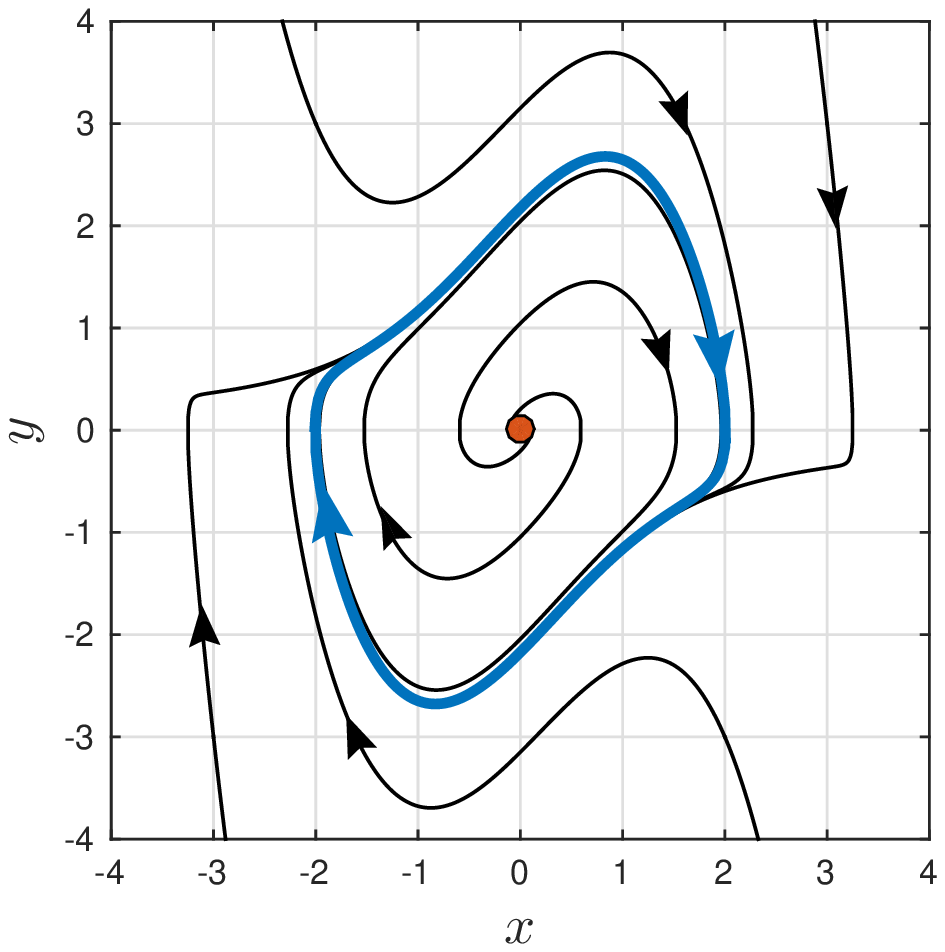}
}
\subfigure[]{\label{F:vdpLimitCycles} 
\includegraphics[width=0.45\textwidth, trim=0cm 3.5cm 5cm 0cm]
{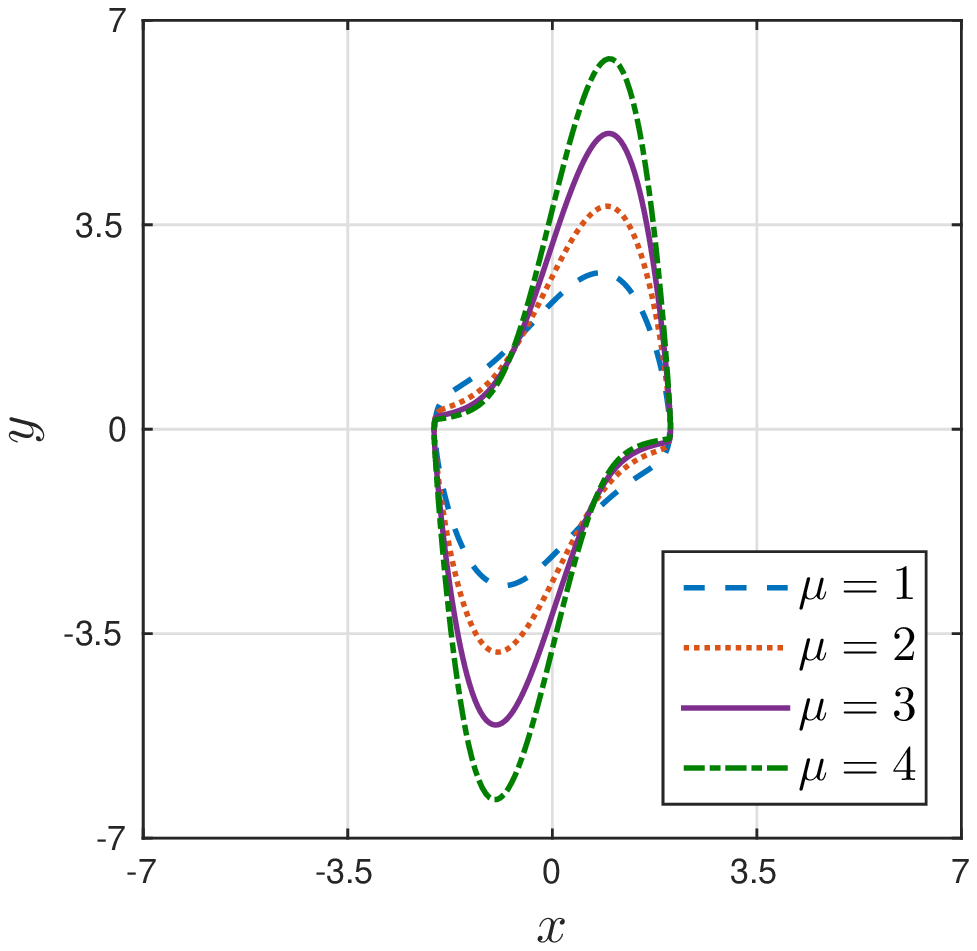}
}
\caption{{(a) Phase portrait of the van der Pol oscillator for $\mu=1$. The 
limit cycle and the unstable fixed point at $(x,y)=(0,0)$ are highlighted. (b) Limit 
cycles for various values of $\mu$.}}
\end{figure}

Throughout this work, we illustrate the methods described using the van der Pol 
oscillator~\cite{Khalil2002}, which can be written as
\begin{equation}
\label{E:VDP_SS}
{\bcolvec{\dot{x}}{\dot{y}}}= \bcolvec{y}{\mu(1-x^2)y-x},
\end{equation}
where a dot denotes $\frac{d}{dt}$. The parameter $\mu>0$ sets the strength of 
the nonlinear damping. There is a limit cycle that attracts all trajectories 
except the unstable fixed point at the origin (Figure~\ref{F:vdpPhasePortrait}), 
and the global invariant set is composed of the limit cycle and the fixed point. 
The system is a standard example of a nonlinear oscillator and has been studied extensively, including with stochastic forcing~\cite{Leung1995,Galan2009,Yuan2013}.
Here we find nearly tight bounds on averages of $x^2+y^2$ both with and without noise.

The rest of this work is organized as follows. Section \ref{SS:SoSBoundsReview} reviews the framework 
of~\cite{Chernyshenko2014a} for bounding deterministic averages and uses it to 
find upper bounds on $\overline{x^2+y^2}$ in the van der Pol system. Section 
\ref{SS:LocalBounds} extends the framework to give attractor-specific bounds, 
which we use to find lower bounds on $\overline{x^2+y^2}$ over the van der Pol 
limit cycle. These lower bounds are larger than zero and thus do not apply to 
the unstable fixed point. The bounding methods for stochastic dynamical systems 
are described in \S\ref{SS:GeneralNoiseFormulation}, and they are specialized to 
the case of small and vanishing noise in 
\S\S\ref{SS:SmallNoise}--\ref{SS:VanishingNoise}. The methods of these sections 
give bounds on $\expect{x^2+y^2}{}$ in the van der Pol example for a range of 
noise amplitudes. Section \ref{S:FurtherComments} discusses the limitations of 
our methods and gives ideas for improvement, and~\S\ref{S:Conclusion} offers 
concluding remarks.

%%%%%%%%%%%%%%%%%%%%%%%% 
%  BOUNDS FOR DETERMINISTIC SYSTEMS
%%%%%%%%%%%%%%%%%%%%%%%%
\section{Global bounds for deterministic systems}
\label{SS:SoSBoundsReview}

Consider an autonomous dynamical system
\begin{equation}
\label{E:DynamicalSystem}
\dot{\x} = \f(\x), 	\quad \x \in \Real^{n},
\end{equation}
in which all trajectories $\x(t)$ remain bounded as $t\to \infty$. We assume nothing special about the structure of $\f$ except that the system is bounded, although eventually we restrict attention to polynomial $\f$. We wish to prove upper and lower bounds on $\overline\phi$, the average of a function 
$\phi(\x)$ over a trajectory $\x(t)$. In applications, $\phi$ can be chosen as a 
quantity of interest for the system. Unless the system has a single attractor that is globally 
asymptotically stable, different trajectories can give different values of 
$\overline\phi$. In this section we seek global bounds on $\overline\phi$, 
meaning they hold for every possible value of $\overline\phi$.

Suppose we wish to prove a constant lower bound $L$ on all possible 
$\overline\phi$,
\begin{equation}
\label{E:LB1}
\overline{\phi} \geq L.
\end{equation}
Central to our method will be a suitably chosen differentiable \textit{storage 
function} $V(\x)$. No matter what $V$ is chosen, $\overline{\f\cdot\del V} = 0$ 
along any bounded trajectory because
\begin{equation}
\overline{\f\cdot\del V} = \overline{\dot V} = 
\lim_{T\to\infty}\frac{1}{T}\left( V\left[\x(T)\right]-V\left[\x(0)\right]\right) 
= 0.
\end{equation}
The desired bound~\eqref{E:LB1} is thus equivalent to the 
inequality\begin{equation}
\label{E:LB2}
\overline{\f\cdot\del{V} + \phi - L} \geq 0
\end{equation}
for any differentiable $V(\x)$. The above time average cannot be evaluated 
without knowing trajectories $\x(t)$, but a sufficient condition for the 
inequality to hold is for it to hold pointwise for all $\x$. Thus, we will have 
proven that $\overline\phi\ge L$ if we can find any differentiable $V(\x)$ such 
that
\begin{equation}
\label{E:LBineq}
%\frac{d V}{d t} + \phi - L \geq 0
\f\cdot\del V + \phi - L \geq 0 \quad \forall\x\in\mathbb R^n.
\end{equation}
An upper bound $\overline\phi\le U$ can be proven in a similar way by reversing 
the inequality sign in the pointwise sufficient condition. The following 
proposition summarizes both conditions.

\begin{proposition}
\label{T:Proposition1}
{
Let $\dot{\x}=\f(\x)$ with $\x\in\mathbb R^n$ be a dynamical system whose 
trajectories are bounded forward in time, let $\phi(\x)$ be a scalar 
function and let $\overline{\phi}$ be its time average defined as in~\eqref{E:TimeAvDef}. If there exist differentiable functions $V_u(\x)$, $V_l(\x)$, and 
constants $U$, $L$ such that
\begin{subequations}
\begin{align}
\label{E:UpperBound}
&& && && \mathcal{D}_u(\x) &:= \f\cdot\del V_u +\phi - U\leq 0  && \forall 
\x\in\Real^{n}, && && && \\
\label{E:LowerBound}
&& && && \mathcal{D}_l(\x) &:=\f\cdot\del V_l + \phi - L \geq 0 && \forall 
\x\in\Real^{n}, && && &&
\end{align}
\end{subequations}
then along any trajectory $\x(t)$,
\begin{equation}
L\leq \overline{\phi}\leq U.
\end{equation}
}
\end{proposition}

%\paragraph{Remark 1}
\begin{remark}
\label{R:Remark1}
In the statement of Proposition~\ref{T:Proposition1} we have assumed that the time average $\overline{\phi}$ exists. Should the limit in~\eqref{E:TimeAvDef} not converge, the Proposition holds if we take the limit superior when computing the upper bound, and the limit inferior when computing the lower bound. That is, for any trajectory $\x(t)$,
\[
L \leq \liminf_{T\to\infty} \int_{0}^{T} \phi[\vec{x}(t)]\,dt \leq  \limsup_{T\to\infty} \int_{0}^{T} \phi[\vec{x}(t)]\,dt \leq U.
\]
See~\cite{Karabacak2011} for a discussion of when infinite-time averages do or do not converge.
\end{remark}

There are two difficulties in applying the above proposition to yield good 
bounds $U$ and $L$. The first is choosing the storage functions $V_u$ and $V_l$. 
The second is checking whether $\mathcal D_u\le0$ and $\mathcal D_l\ge0$ for 
candidate storage functions and bound values. These difficulties can be 
prohibitive in general, but the task is greatly simplified when $\f$ and $\phi$ 
are polynomials of the state variables $x_1, \hdots ,x_n$.
%$x_i$, $i\in\set{1,...,n}$. 

In the rest of this work we assume that $\f$ and $\phi$ are polynomials. If the 
chosen $V_u$ and $V_l$ are also polynomials, then so are $\mathcal{D}_u$ and 
$\mathcal{D}_l$. Checking the sufficient conditions~\eqref{E:UpperBound} 
and~\eqref{E:LowerBound} thus amounts to verifying the non-negativity of 
polynomial expressions. While this is an NP-hard problem, the computational 
complexity can be significantly reduced by replacing the conditions 
$\mathcal{D}_l(\x)\geq 0$ and  $-\mathcal{D}_u(\x)\geq 0$ with the stronger 
conditions that $\mathcal{D}_l$ and $-\mathcal{D}_u$ belong to the set $\Sigma$ 
of polynomials that are sums-of-squares (SoS). A polynomial $P(\x)$ is said to 
be SoS if it can be expressed as the sum of squares of some other polynomials---that is, if there exist polynomials $\set{p_i(\x)}_{i=1}^M$ such that
\begin{equation}
P(\x) = \sum_{i=1}^{M} p_i(\x)^2.
\end{equation}

To prove the bounds $U$ and $L$ it suffices to find polynomials 
$V_u$ and $V_l$ such that $\mathcal{D}_l\in\Sigma$ and 
$-\mathcal{D}_u\in\Sigma$. The best bounds that can be proven in this framework are
\begin{align}
\label{E:OptBound}
&\begin{gathered}
\min_{V_u} \quad U		\\
\text{s.t.} \quad -\left(\f\cdot\del V_u + \phi - U\right)  \in \Sigma,
\end{gathered}
&
&\begin{gathered}
\max_{V_l} \quad L		\\
\text{s.t.} \quad \f\cdot\del V_l + \phi - L \in \Sigma.
\end{gathered}
\end{align}
For the storage functions $V_u$ and $V_l$, we must specify polynomial ansatze with undetermined 
coefficients. The decision variables in the upper bound optimization, for 
instance, are $U$ and the coefficients of $V_u$.

\paragraph{Computational methods.} 
Optimization problems with SoS constraints such as~\eqref{E:OptBound} can be 
solved numerically using the methods of SoS programming.
The main idea behind SoS programming is that every polynomial can be represented 
as a symmetric matrix (after defining a suitable polynomial basis set), and this 
matrix can be positive semidefinite if and only if the polynomial admits a SoS 
decomposition~\cite{Parrilo2000,Parrilo2003,Blekherman2012}. Constraints involving SoS 
polynomials, including additional equality and inequality constraints, can be 
posed as equality and inequality conditions on symmetric matrices. An 
optimization problem with constraints of this type is known as a semidefinite 
program (SDP). A number of efficient computer solvers are available for SDPs 
(e.g.~\cite{Sturm1999,Toh1999,Tutuncu2003,Fujisawa2008,Andersen2009}), and the 
software packages YALMIP~\cite{Lofberg2004} or 
SOSTOOLS~\cite{Papachristodoulou2013} can assist in formulating SoS constraints 
as SDPs. More details on convex optimization and the link between SoS 
polynomials and SDPs can be found in~\cite{Boyd2004,Parrilo2000,Lofberg2009}, 
while examples using SoS programming to study dynamical systems can be found 
in~\cite{Parrilo2000,Papachristodoulou2002,Papachristodoulou2005,Tan2006a,
Goulart2012,Chernyshenko2013,Henrion2014,Valmorbida2014,Huang2015}.

For our numerical implementation, we used the SoS module of 
YALMIP~\cite{Lofberg2009} to transform SoS optimization problems into SDPs. To 
solve the SDPs, we used the multiple-precision solver 
SDPA-GMP~\cite{Fujisawa2008} for the following reasons.
First, the SDPs we solved were badly conditioned even for modest polynomial 
degrees, and none of the standard double-precision solvers we tested converged 
reliably. This issue could be resolved by carefully rescaling the system, but 
there is no systematic rule to determine a suitable rescaling. Second, working 
in multiple-precision offers a simple check on results: increase the precision, 
and confirm that the bounds $U$ or $L$ change very little. We have done this for 
all bounds reported here. Appendix~\ref{S:RigorousBounds} describes additional 
checks on numerical results, including how SDP computations could be part of 
fully rigorous computer-assisted proofs.

%%%%%%%%%%%%%%%%%%%%%%%%
\paragraph{Example\,{\rm{:}} deterministic upper bounds for the van der Pol 
limit cycle.}
\label{SS:ExampleDeterministicVDP_UB}

\begin{figure}%[t!]
\centering
\includegraphics[width=0.60\textwidth, trim=0cm 0cm 0cm 1cm]{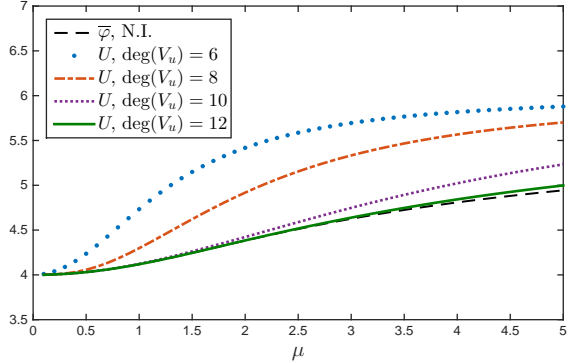}
\caption{{Optimal upper bounds on $\overline{\phi}=\overline{x^2+y^2}$ for 
the van der Pol oscillator computed with the upper bound problem 
of~\eqref{E:OptBound} for different degrees of $V_u$. The time average 
$\overline{\phi}$ obtained by numerical integration (N.I.) 
of the system is also shown.}}
\label{F:UBvdpSDPAGMP} 
\end{figure}

To illustrate the application of Proposition~\ref{T:Proposition1}, we
have computed upper bounds $\overline\phi=\overline{x^2+y^2}\le U$ in the van 
der Pol system~\eqref{E:VDP_SS} for various $0.1\leq \mu\leq 5$ and four 
different polynomial degrees of $V_u$. Figure~\ref{F:UBvdpSDPAGMP} shows the 
resulting $U$, along with estimates of $\overline\phi$ obtained by 
integrating over the limit cycle using the fourth-order Runge-Kutta method with 
fixed time steps. As expected, increasing the degree of the storage function 
$V_u$ gives a better (smaller) upper bound $U$. For a given degree, the bounds 
become less tight as $\mu$ is raised because the shape of the limit cycle is 
more complicated (cf. Figure~\ref{F:vdpLimitCycles}).

The trivial lower bound $0\le\overline{x^2+y^2}$ is the best \textit{global} 
lower bound possible for the van der Pol system. It cannot be improved because 
it is saturated by the trajectory staying at the unstable fixed point 
$(x,y)=(0,0)$. However, if one is interested only in the trajectories that tend to the stable limit cycle, then the lower bound of zero 
is not tight. 
Raising the lower bound requires a result that does not apply globally, but instead applies only to a subset of all possible trajectories.

%%%%%%%%%%%%%%%%%%%%%%%% 
%  BOUNDS ON LOCAL ATTRACTORS
%%%%%%%%%%%%%%%%%%%%%%%%
\section{Bounds on local attractors}
\label{SS:LocalBounds}

Often, one is interested in the average of a function $\phi(\x)$ over only a 
subset of all possible trajectories, such as those tending to a particular 
attractor. If a bound on $\overline\phi$ over such trajectories is to be tight, 
it should not apply to trajectories that tend to a different local attractor, 
nor should it apply to trajectories on unstable invariant structures that are 
not part of the attractor of interest (such as unstable fixed points, unstable 
limit cycles, or basin boundaries). However, global bounds of the kind derived 
using Proposition~\ref{T:Proposition1} are unlikely to be tight for a particular 
attractor since they must obey
\begin{equation}
L \le \inf_{\x(t)}\overline{\phi(\x(t))} \le \sup_{\x(t)}\overline{\phi(\x(t))} 
\le U,
\end{equation}
where $\x(t)$ is any trajectory of the dynamical system. 

Suppose we wish to bound possible values of $\overline\phi$, not for all 
trajectories but only for trajectories that eventually enter and remain inside a 
given absorbing domain $\mathcal T$. To compute $\overline\phi$ over any such 
trajectory, it suffices to begin averaging after the trajectory has permanently 
entered $\mathcal T$, so dynamics outside of $\mathcal T$ can be ignored. This 
suggests the following modification of Proposition \ref{T:Proposition1}, where 
the inequality conditions are imposed only on $\mathcal T$ and, consequently, 
the resulting bounds are proven only for trajectories that permanently enter 
$\mathcal T$. The comments made in Remark~\ref{R:Remark1} still apply.

\begin{proposition}
\label{T:Proposition2}
Let $\dot{\x}=\f(\x)$ be a dynamical system with $\x\in\mathbb R^n$, let 
$\phi(\x)$ be a scalar function, let $\overline{\phi}$ be its time average defined as in~\eqref{E:TimeAvDef}, and let $\mathcal T$ be an absorbing domain. If 
there exist differentiable functions $V_u(\x)$, $V_l(\x)$, and constants $U$, $L$ 
such that
\begin{subequations}
\begin{align}
\label{E:UpperBound_Sproc}
&& && && \mathcal{D}_u(\x) &:= \f\cdot\del V_u +\phi - U\leq 0  && \forall 
\x\in\mathcal{T}, && && && \\
\label{E:LowerBound_Sproc}
&& && && \mathcal{D}_l(\x) &:=\f\cdot\del V_l + \phi - L \geq 0 && \forall 
\x\in\mathcal{T}, && && &&
\end{align}
\end{subequations}
then along any trajectory $\x(t)$ that permanently enters $\mathcal T$,
\begin{equation}
L\leq \overline{\phi}\leq U.
\end{equation}
\end{proposition}

Applying Proposition \ref{T:Proposition2} to a particular system requires 
specifying an appropriate absorbing domain. For instance, if one is interested 
in a local attractor $\mathcal A$ with basin $\mathcal B$, it suffices to choose 
any $\mathcal T$ such that $\mathcal A\subset\mathcal T\subset\mathcal B$. 
Bounds proven on $\mathcal T$ then apply to all trajectories that approach the 
attractor or are part of the attractor, but they need not apply to trajectories 
outside $\mathcal B$. 

Finding a mathematical description of $\mathcal T$ for a given attractor is 
difficult in general, but this too can be done using SoS programming 
\cite{Tan2006a,Henrion2014}. An example in which it is \emph{not} difficult to 
choose $\mathcal T$ is when the trajectory to be excluded from the bound is a 
repelling fixed point. In this case it suffices to choose $\mathcal T=\mathbb 
R^n\smallsetminus\mathcal U$ for any small enough set $\mathcal U$ around that 
point.

If the absorbing domain can be specified as a semi-algebraic set---that is, 
defined by a set of polynomial inequalities and equalities---the conditions of 
Proposition~\ref{T:Proposition2} can be checked by SoS programming using the 
\textit{generalized \Sprocedure}~\cite[Lemma 2.1]{Tan2006}. For instance, 
suppose $\mathcal{T} = \set{\x \,|\, g(\x)\geq 0}$ for some polynomial $g$; 
proving a lower bound for trajectories entering $\mathcal T$ calls for SoS 
conditions that imply $\mathcal D_l(\x)\geq 0$ when $g(\x)\geq 0$.  
For this to be true, it suffices that there exist $s(\x)\ge0$ such that 
$\mathcal{D}_l(\x) - s(\x)g(\x)\geq 0$. Strengthening these non-negativity 
constraints to SoS constraints and making similar arguments for the upper bound 
gives the following two SoS programs:
\begin{align}
\label{E:OptBoundSproc}
&\begin{gathered}
\min_{V_u,s} \quad U		\\
\text{s.t.} \quad -\left(\f\cdot\del V_u + \phi - U\right) - s\,g \in \Sigma \\
s \in \Sigma,
\end{gathered}
&
&\begin{gathered}
\max_{V_l,s} \quad L		\\
\text{s.t.} \quad \f\cdot\del V_l + \phi - L - s\,g \in \Sigma \\
s \in \Sigma,
\end{gathered}
\end{align}
where polynomial ansatze are specified for $V_u$, $V_l$, and $s$, and their free coefficients are the decision variables.
Using similar ideas, the \Sprocedure~can be generalized to semi-algebraic 
$\mathcal T$ defined by multiple polynomial inequalities and equality 
constraints~\cite{Tan2006}.

%%%%%%%%%%%%
\paragraph{Example\,{\rm{:}} deterministic lower bounds for the van der Pol 
limit cycle.}
\begin{figure}[b!]
\centering
\subfigure[]{\label{F:LBvdpR05} 
\includegraphics[width=0.47\textwidth, trim=0cm 0.5cm 1cm 1cm]
{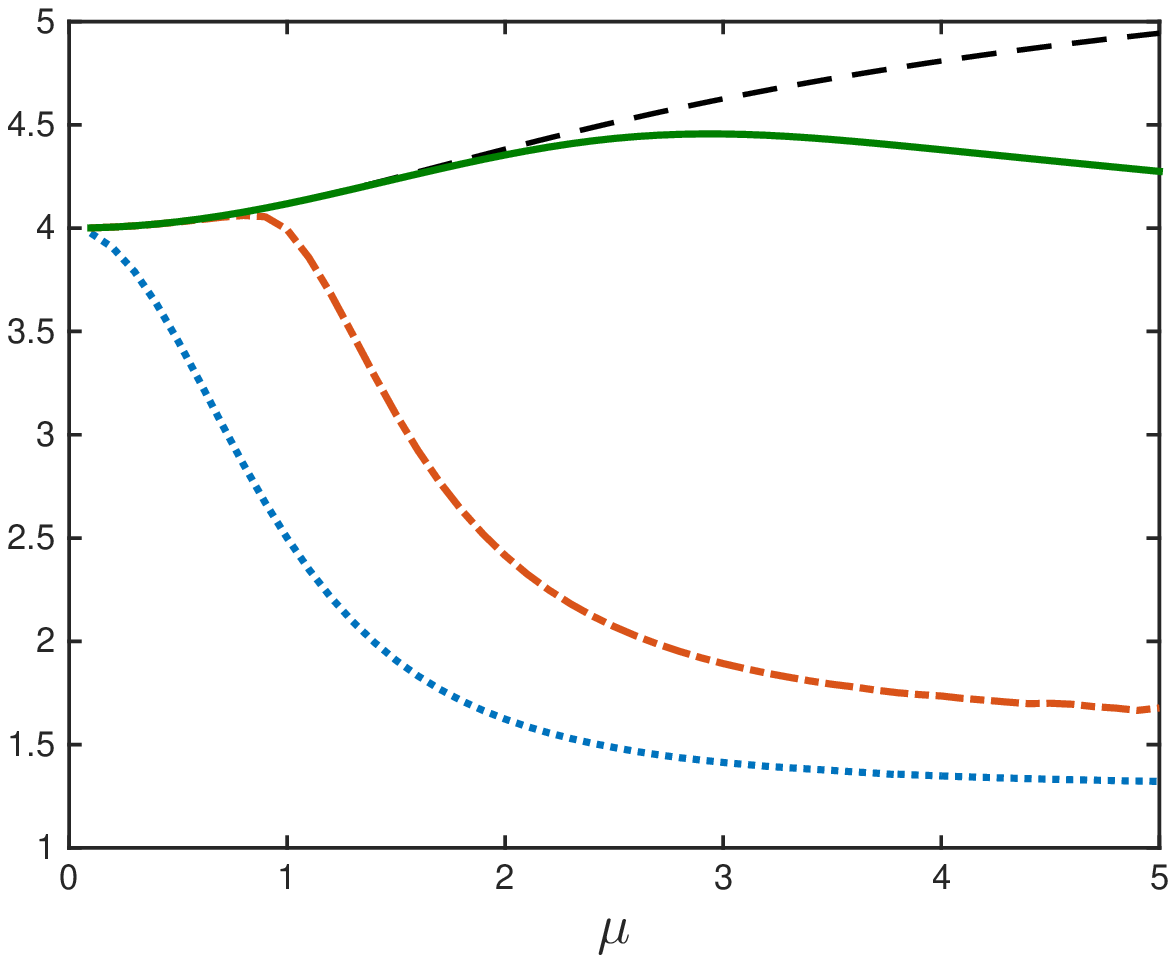} }
\hfill
\subfigure[]{\label{F:LBvdpR1} 
\includegraphics[width=0.47\textwidth, trim=1cm 0.5cm 0cm 1cm]
{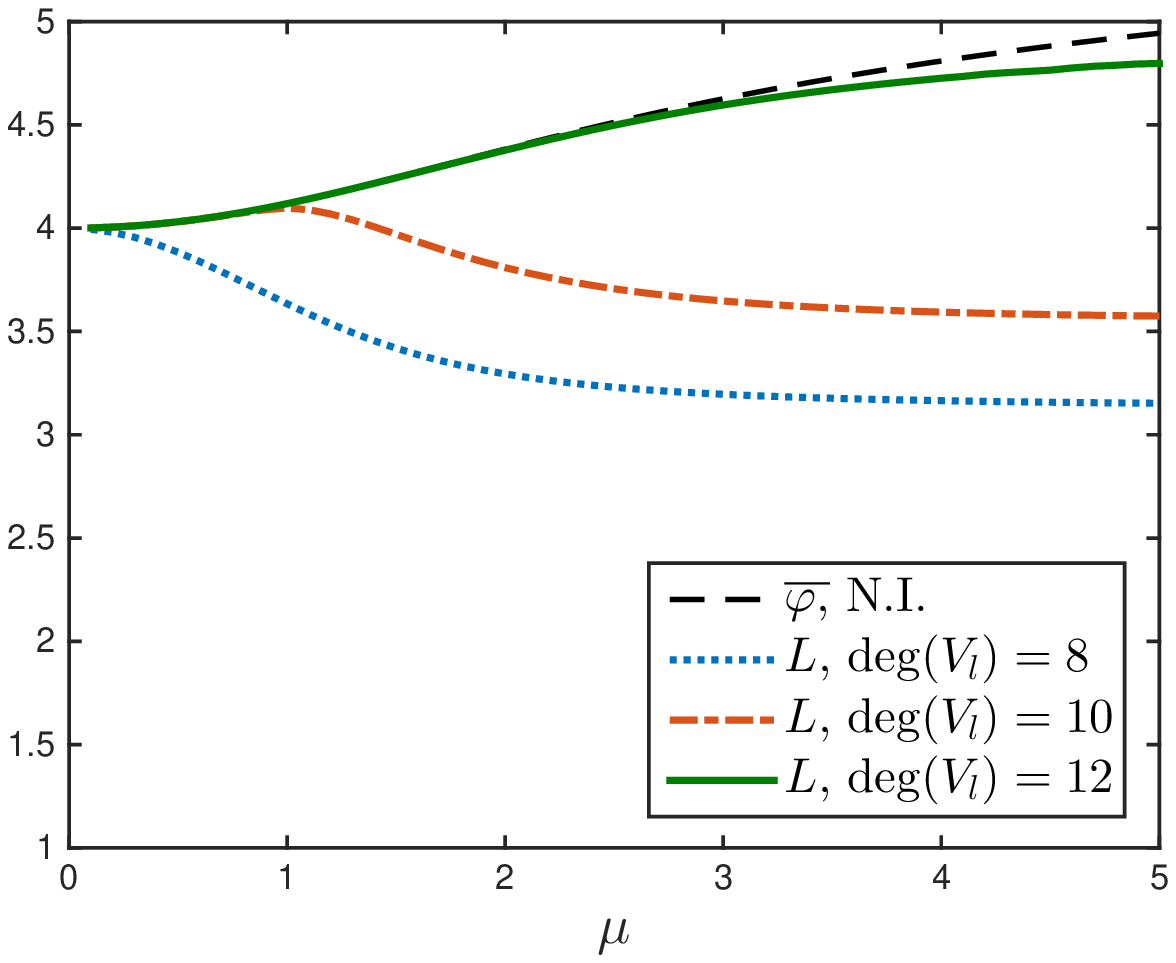} }
\caption{{Optimal lower bounds on $\overline{\phi}=\overline{x^2+y^2}$ over 
the limit cycle of the  van der Pol system for different degrees of $V_l$ using 
the absorbing domains (a) 
$\mathcal{T}_{0.5} = \set{(x,y) \, |\, x^2+y^2-0.25 \geq 0 }$ and (b) 
$\mathcal{T}_{1} = \set{(x,y) \, |\, x^2+y^2-1 \geq 0 }$. The time average 
$\overline{\phi}$ from numerical integration (N.I.) is also shown.}}
\label{F:LBvdp}
\end{figure}
Let us revisit the lower bound on $\overline{x^2+y^2}$ for the van der Pol 
oscillator. Suppose we want a bound that does not apply to the unstable fixed 
point at the origin but applies to the other trajectories, all of which approach 
the limit cycle. If such a lower bound is perfectly tight, it will equal the 
value of $\overline\phi$ on the limit cycle. To apply 
Proposition~\ref{T:Proposition2}, we must specify an absorbing domain $\mathcal 
T$ that contains the entire limit cycle but omits the origin. As described in 
Appendix~\ref{S:VDPAbsorbingDomainProof}, the set $\mathcal{T}_r = \set{(x,y) \, 
|\, g(x,y)=x^2+y^2-r^2 \geq 0 }$ is just such an absorbing domain for any 
$r\leq1$.

Figure~\ref{F:LBvdp} shows the lower bounds obtained by solving the lower bound 
program of~\eqref{E:OptBoundSproc} using the absorbing domains 
$\mathcal{T}_{0.5}$ and $\mathcal{T}_{1}$. At a given polynomial degree, using 
the smaller absorbing domain $\mathcal{T}_{1}$ gives a better bound. Results 
would likely be improved further by using $\mathcal T$ that closely approximate 
the attractor. This suggests a two-step procedure: using SoS techniques to find 
the smallest possible $\mathcal T$ around the attractor of 
interest~\cite{Tan2006a,Henrion2014}, and then solving the SoS bounding 
programs~\eqref{E:OptBoundSproc}.

%%%%%%%%%%%%%%%%%%%%%%%%
% BOUNDS FOR STOCHASTIC SYSTEMS
%%%%%%%%%%%%%%%%%%%%%%%%
\section{Bounds for stochastic systems}
\label{SS:GeneralNoiseFormulation}

%When dealing with physical systems, it is usual to account for the effect of external disturbances or unmodeled physical mechanisms by adding noise to the deterministic dynamics~\cite{Doering1990}. For this reason, we will now turn our attention to the problem of determining bounds on expectations for stochastic systems.
%
%Specifically, 
Consider a stochastic dynamical system
\begin{equation}
\label{E:StochasticDynamicalSystem}
\dot{\x}=\f(\x) + \sqrt{2\e}\,\boldsymbol{\sigma} \boldsymbol{\xi},	\quad 
\x\in \Real^{n},\, \boldsymbol\xi \in \Real^{m},
\end{equation}
in which the random trajectories $\x(t)\in\Real^n$ are bounded almost surely as 
$t\to\infty$. The standard vector Wiener process 
$\boldsymbol{\xi}(t) \in \Real^{m}$ is 
pre-multiplied by the matrix $\boldsymbol{\sigma} \in \Real^{n\times m}$ that 
describes the relative effect of each noise component ${\xi}_i$ on each state 
variable $x_j$, and the overall noise strength is scaled by $\sqrt{2\e}$. 

%Interest in stochastic systems ofthe form~\eqref{E:StochasticDynamicalSystem} generally stems from the need to model the effect of external disturbances. 
We assume for simplicity that the noise is additive, meaning that 
$\boldsymbol\sigma$ does not vary with $\x$, but the present analysis extends 
with only minor modifications to multiplicative noise in either the It\^o or Stratonovich 
interpretation.
We also assume that the system has reached statistical equilibrium, in which case the 
probability density of its trajectories, $\rho(\x)\ge0$, decays at infinity and 
satisfies the stationary Fokker--Planck equation
\begin{equation}
\label{E:SteadyFP}
\del\cdot( \e \vec{D}\del\rho -  \f  \rho) = 0, \qquad 
\int_{\Real^n}\rho(\x)d\x=1, \qquad 
\vec{D}:=\boldsymbol{\sigma}\boldsymbol{\sigma}^T.
\end{equation} 

Suppose we wish to prove a constant lower bound $\expect{\phi}{\e}\geq L$, where 
$\expect{\phi}{\e}$ is the stationary expectation of $\phi(\x)$  defined as 
in~\eqref{E:ExpectDef}. We assume that such a stationary average exists, which is true for all $\phi(\x)$ that don't grow too fast as $\vert \x \vert\to\infty$ (precise statements can be found in~\cite{Meyn1993}). The subscript on $\expect{\phi}{\e}$ indicates its 
dependence on the noise strength $\e$. Since $\expect{L}{\e}=L$, this is 
equivalent to proving
\begin{equation}
\label{E:StochBound1}
\expect{\phi - L}{\e} \geq 0.
\end{equation}
As in the deterministic case, our method of proof relies on a suitably chosen 
differentiable storage function $V(\x)$. For any $V$ that does not grow too 
quickly as $|\x|\to\infty$, the stationary expectation $\expect{\e\del 
\cdot(\vec{D} \del V) + \f\cdot\del V}{\e}$ is zero because
\begin{equation}
\label{E:zero}
\begin{aligned}
\expect{\e\del \cdot(\vec{D} \del V) + \f\cdot\del V}{\e} 
	&= \int_{\Real^n}\rho\left[\e \del \cdot(\vec{D} \del V) + \f\cdot\del V 
\right]d\x \\
	&= \int_{\Real^n}V\,\del\cdot( \e \vec{D}\del\rho -  \f  \rho)d\x \\
	&= 0. 
\end{aligned}
\end{equation}
The third line above follows from the stationary Fokker-Planck 
equation~\eqref{E:SteadyFP}. The second line follows from integration by parts, 
assuming that the boundary terms vanish---that is, assuming 
\begin{equation}
\label{E:boundaryInt}
\lim_{R\to\infty}\int\limits_{\vert\x\vert=R} \left( \e\, \rho\, \vec{D}\del V  
- \e \,V\, \vec{D}\del \rho + \rho \,V\,\f \right)\cdot \boldsymbol{\nu} \,dS = 
0,
\end{equation}
where $\boldsymbol{\nu}(\x)$ is the outwards unit normal to the sphere 
$\vert\x\vert=R$, and $dS$ is the surface element. 
The above condition holds in many cases, including any case where $V$ is 
polynomial and $\rho$ decays exponentially at infinity. When this condition 
holds, then so does the equality~\eqref{E:zero}, hence the 
inequality~\eqref{E:StochBound1} is equivalent to 
\begin{equation}
\label{E:LB_IntIneq}
\expect{ \e\del \cdot(\vec{D} \del V) + \f\cdot\del V + \phi - L}{\e} \ge 0.
\end{equation}
The above expectation cannot be evaluated without knowing $\rho$, but it is 
sufficient for the inequality to hold pointwise for all $\x$. The lower bound 
$\expect{\phi}{\e}\geq L$ is therefore proven if we can find any differentiable 
$V(\x)$ satisfying the boundary integral condition~\eqref{E:boundaryInt} and the 
pointwise inequality
\begin{equation}
\e\del \cdot(\vec{D} \del V) + \f\cdot\del V + \phi - L \geq 0 \quad 
\forall\x\in\mathbb R^n.
\end{equation}

The same argument with a reversed inequality sign gives a sufficient condition 
for an upper bound $U$ on $\expect{\phi}{\e}$. We summarize these results in the 
following Proposition, which is the stochastic analog of 
Proposition~\ref{T:Proposition1}.
%
%\vskip 10pt
\begin{proposition}
\label{T:Theorem3}
Let $\dot{\x} = \f(\x)+\sqrt{2\e}\,\boldsymbol{\sigma} \boldsymbol{\xi}$ with 
$\x\in\Real^{n}$, $\boldsymbol{\xi}\in\Real^{m}$, $\boldsymbol{\sigma}\in\Real^{n\times m}$ be 
a stochastic dynamical system in a statistically stationary state with 
probability distribution $\rho(\x)$. If there exist differentiable functions 
$V_u$, $V_l$ and constants $U$, $L$ such that
\begin{subequations}
\label{E:StochasticInequalities}
\begin{align}
\label{E:V_Noise_INEQ_UB}
&& && && \e \del \cdot(\vec{D} \del V_u) + \f\cdot\del V_u + \phi - U  &\leq 0 
&&\forall \x \in \Real^n,	&& && &&\\
\label{E:V_Noise_INEQ_LB}
&& && && \e \del \cdot(\vec{D} \del V_l) + \f\cdot\del V_l + \phi - L  &\geq 0 
&&\forall \x \in \Real^n, && && &&
\end{align}
\end{subequations}
where $\vec{D}= \boldsymbol{\sigma}\boldsymbol{\sigma}^T$, and if $V_u$, $V_l$ 
grow slowly enough at infinity to each satisfy~\eqref{E:boundaryInt}, then
\begin{equation}
L\leq \expect{\phi}{\e}\leq U.
\end{equation}
\end{proposition}

The inequality conditions~\eqref{E:V_Noise_INEQ_UB} 
and~\eqref{E:V_Noise_INEQ_LB} were derived by a different method 
in~\cite{Chernyshenko2014a} for the case where $\vec D$ is the identity matrix 
(hence $\del \cdot\vec{D} \del=\del^2$). The conditions differ from their 
deterministic counterparts~\eqref{E:UpperBound} and~\eqref{E:LowerBound}, 
respectively, only in the addition of the diffusive terms $\e \del \cdot(\vec{D} 
\del V_u)$ and $\e \del \cdot(\vec{D} \del V_l)$. Similar results for non-negative $\phi$ were also stated in~\cite{Glynn2008}. 

Like the deterministic bounds 
of~\S\S\ref{SS:SoSBoundsReview}--\ref{SS:LocalBounds}, the stochastic bounds and 
storage functions in Proposition~\ref{T:Theorem3} can be found numerically for a 
given $\e$ if $\f$ and $\phi$ are polynomials. Letting $V_u$ and $V_l$ be 
polynomials also, we replace~\eqref{E:V_Noise_INEQ_UB} 
and~\eqref{E:V_Noise_INEQ_LB} with stronger SoS constraints to obtain the SoS 
programs
\begin{align}
\label{E:OptBoundFixedNoise}
&\begin{gathered}
\min_{V_u} \quad U		\\
\text{s.t.} \; -\left[\e \del \cdot(\vec{D} \del V_u) + \f\cdot\del V_u + \phi - 
U \right]\in \Sigma,
\end{gathered}
&
&\begin{gathered}
\max_{V_l} \quad L		\\
\text{s.t.} \; \e \del \cdot(\vec{D} \del V_l) + \f\cdot\del V_l + \phi - L\in 
\Sigma.
\end{gathered}
\end{align}
Since Proposition \ref{T:Theorem3} relies on the boundary integral 
condition~\eqref{E:boundaryInt}, the $V_u$ and $V_l$ constructed by the above 
SoS programs must satisfy this condition in order for $U$ and $L$ to be proven 
bounds. One way to guarantee this is to prove that $\rho$ decays exponentially 
as $|\x|\to\infty$.

%%%%%%%%%%%%%%
\paragraph{Example\,{\rm{:}} stochastic bounds for the van der Pol oscillator.}
\label{SS:ExampleNoisyVDPnoLog}

To illustrate the application of Proposition \ref{T:Theorem3}, we have bounded 
the stationary expectation $\expect{\phi}{\e}=\expect{x^2+y^2}{\e}$ for the 
stochastic van der Pol oscillator
\begin{equation}
\label{E:VDP_Noise}
%\frac{d}{dt}\pcolvec{x}{y}= \pcolvec{y}{\mu(1-x^2)y-x}.
{\bcolvec{\dot{x}}{\dot{y}}}= \bcolvec{y}{\mu(1-x^2)y-x} + 
\sqrt{2\e}\begin{bmatrix}
0 \\ 1
\end{bmatrix} \xi
%\ddot{x} - \mu(1-x^2)\dot{x} + x = \sqrt{2\e}\, \xi,
\end{equation}
with $\mu=1$. 
%The same system was considered e.g. in~\cite{Leung1995} as an example of stochastic system with limit-cycle dynamics. 
For the $\boldsymbol\sigma$ matrix chosen above, the Wiener 
processes $\xi$ acts on the $y$ component alone. This corresponds to a stochastic physical force, as seen by rewriting~\eqref{E:VDP_Noise} as a second-order equation for the position $x(t)$,
\begin{equation}
\ddot{x} - \mu(1-x^2)\dot{x} + x = \sqrt{2\e} \,\xi.
\end{equation}
\begin{figure}%[b!]
\centering
\subfigure[]{\label{F:UB_Noisy_vdp} 
\includegraphics[width=0.47\textwidth, trim=0cm 3.5cm 5cm 0.5cm]
{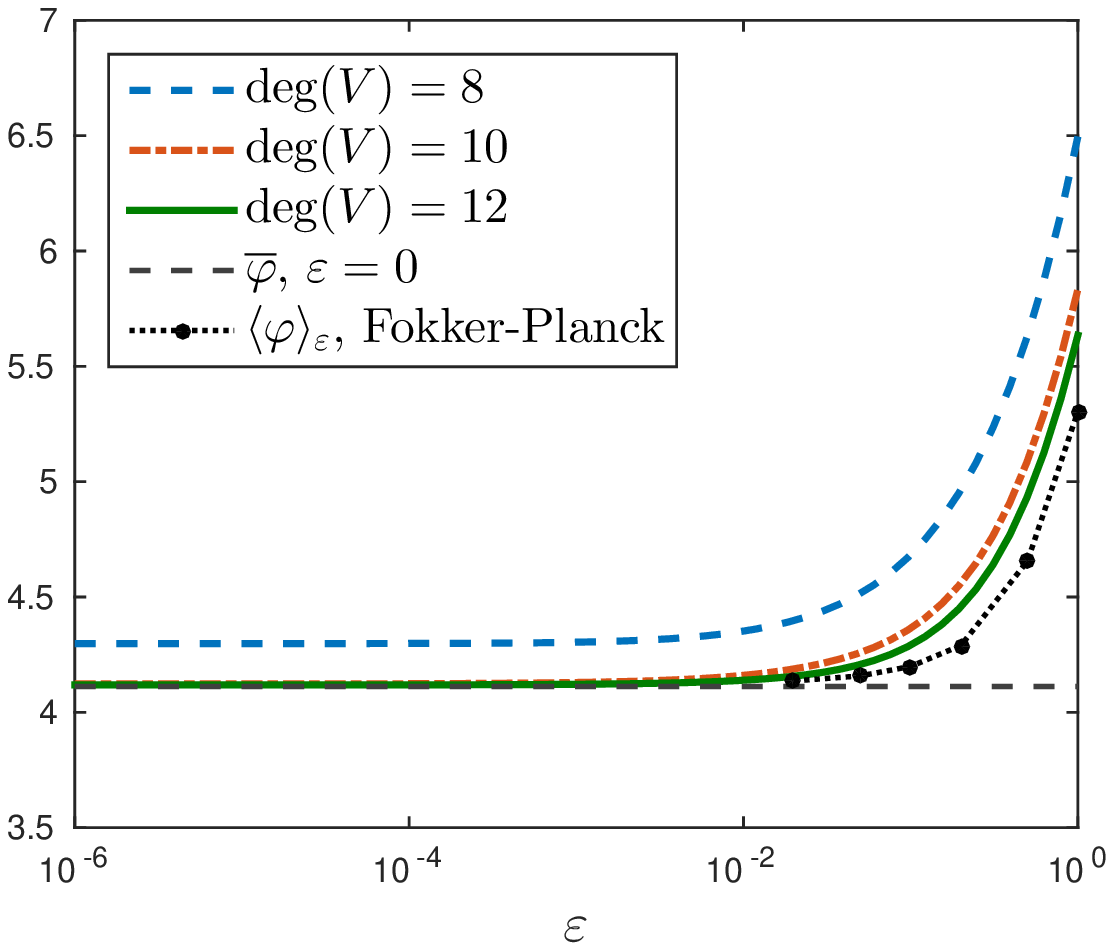} }
\subfigure[]{\label{F:LB_Noisy_vdp} 
\includegraphics[width=0.47\textwidth, trim=0cm 3.5cm 5cm 0.5cm]
{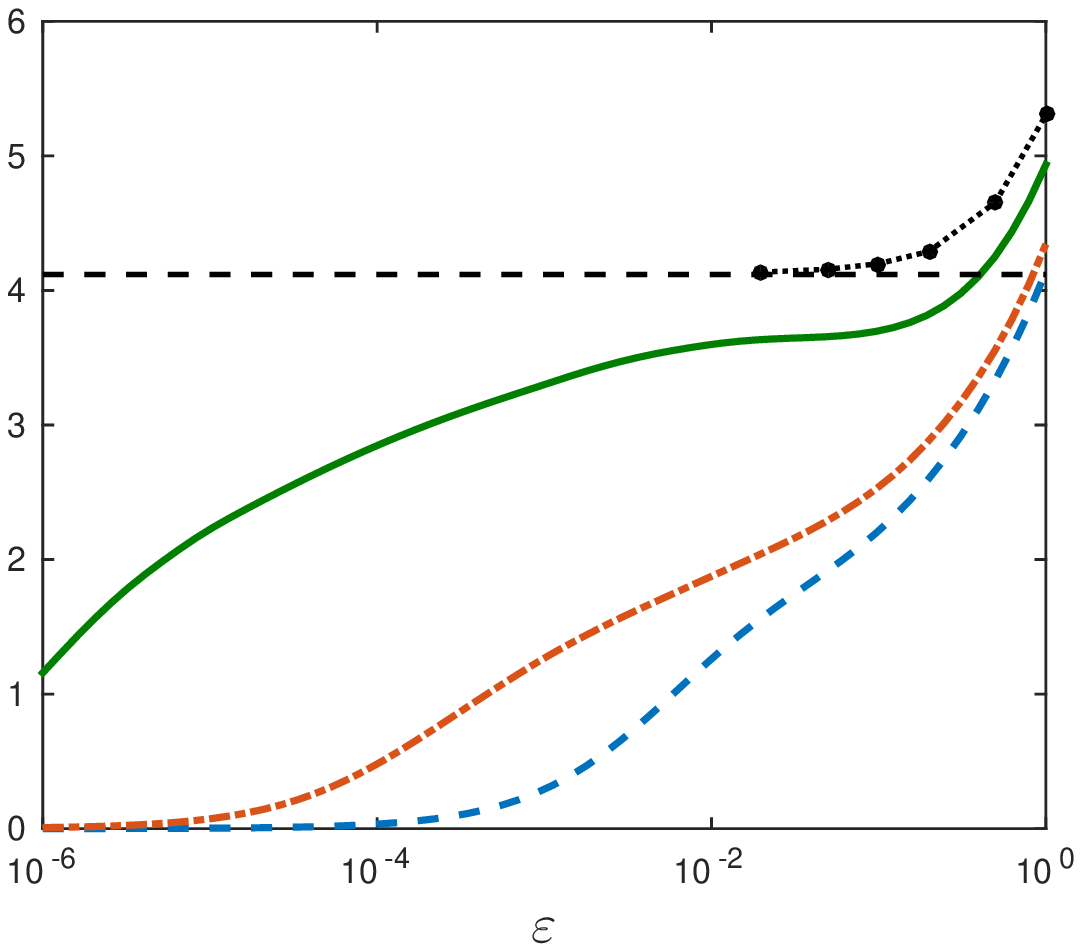} }
\caption{(a) Upper bounds and (b) lower bounds on 
$\expect{\phi}{\e} = \expect{x^2+y^2}{\e}$ for the stochastic van der Pol 
oscillator for $\mu=1$ as a function of the noise strength $\e$. 
The deterministic ($\e=0$) 
average $\overline{\phi}=\overline{x^2+y^2}\approx 4.118$ and the 
expectation $\expect{\phi}{\e}$ computed from the 
Fokker-Planck equation for selected values of $\e$ are also shown. Improved 
lower bounds appear in Figure \ref{F:NoisyVDP_Log}.}
\label{F:NoisyVDP}
\end{figure}
We have computed bounds using the SoS programs~\eqref{E:OptBoundFixedNoise} for 
noise amplitudes $10^{-6} \leq \e \leq 1$ and three different polynomial degrees 
of $V_u$, $V_l$. Figure~\ref{F:NoisyVDP} shows the optimal bounds, along with 
several values of $\expect\phi\e$ computed by numerically solving the stationary 
Fokker-Planck equation~\eqref{E:SteadyFP}. Our solution of the Fokker-Planck 
equation employed finite differences with operator splitting as in 
\cite{Chang1970}, and the system was evolved to steady state by implicit Euler 
time stepping. The limiting expectation $\lim_{\e\to 0}\expect{\phi}{\e}$ equals 
the deterministic average $\overline{\phi}$ on the limit cycle, reflecting the 
fact that the van der Pol system is stochastically stable.

The upper bounds $U$ in Figure~\ref{F:UB_Noisy_vdp} are nearly tight at all noise 
strengths for $V_u$ of degree 10 and higher. They become less tight as $\e$ 
increases, but they still correctly capture the increase in 
$\expect{x^2+y^2}{\e}$ that occurs when stochastic forcing ``smears out'' the 
van der Pol limit cycle.

The lower bounds $L$ in Figure~\ref{F:LB_Noisy_vdp} are good for strong noise 
but not for weak noise. In fact, $L\to0$ as $\e\to0$ for any fixed polynomial 
degree of $V_l$. To see that this is inevitable, observe that the diffusive term 
$\e \del \cdot(\vec{D} \del V_l)$ in \eqref{E:V_Noise_INEQ_LB} vanishes as 
$\e\to0$ if $\del \cdot(\vec{D} \del V_l)$ is bounded. 
As $\e\to0$, the stochastic SoS programs~\eqref{E:OptBoundFixedNoise} reduce to 
the deterministic SoS programs~\eqref{E:OptBound}, so the stochastic bounds can 
be no tighter than the global deterministic bounds. 

In the van der Pol example, upper bounds are not affected by this phenomenon 
since the deterministic global upper bound on $\overline{\phi}$ of 
\S\ref{SS:ExampleDeterministicVDP_UB} is sharp and $\expect\phi\e \to 
\overline{\phi}$ as $\e\to0$. However, a tight lower bound on $\expect\phi\e$ at small 
$\e$ is impeded by the unstable deterministic trajectory at the 
origin, which saturates the deterministic global lower bound $\overline\phi \geq 
0$. In the deterministic case, a tight local lower bound was achieved in 
\S\ref{SS:LocalBounds} by ignoring part of phase space, but we cannot do so in 
the stochastic case since $\rho>0$ everywhere. Instead, the diffusive term 
$\e\del \cdot(\vec{D} \del V_l)$ in the SoS constraint 
of~\eqref{E:OptBoundFixedNoise} must remain $O(1)$ near the origin as $\e\to0$. 
Clearly, this cannot be accomplished with polynomial $V_l$. The next 
section illustrates how a non-polynomial $V_l$ can  be used instead to improve the 
lower bounds on $\expect\phi\e$ at small $\e$.

%%%%%%%%%%%%%%%%%%%%%%%% 
%  BOUNDS WITH WEAK NOISE
%%%%%%%%%%%%%%%%%%%%%%%%
\section{Bounds for stochastic systems with weak noise}
\label{SS:SmallNoise}

An unstable trajectory of the deterministic system $\dot{\x}=\f(\x)$ can 
interfere not only with tight bounds on $\overline\phi$ for a local attractor 
but also with tight bounds on $\expect\phi\e$ in the corresponding stochastic 
system, at least when the noise is weak. This was illustrated in the previous 
section, where lower bounds on $\expect{x^2+y^2}\e$ in the van der Pol system 
were not tight: they approached zero as $\e\to0$ with polynomial $V_l$ of fixed 
degree. In this section we derive lower bounds on $\expect{x^2+y^2}\e$ that 
remain tight when $\e$ is small. We do this using a non-polynomial $V_l$ that 
depends explicitly on $\e$.

The method described below applies to stochastic upper or lower bounds, provided 
that the deterministic trajectory interfering with the bounds is a 
\emph{repelling fixed point}. For concreteness, suppose we are trying to prove a 
tight lower bound on an expectation $\expect\phi\e$ that is strictly positive, 
as in the van der Pol example. Suppose also that the global lower bound on 
$\overline\phi$ in the corresponding deterministic system, $\dot{\x}=\f(\x)$, is 
zero and is saturated by the repelling fixed point $\x=\vec{0}$. For the lower 
bound on $\expect\phi\e$ to remain larger than the deterministic lower bound as 
$\e\to 0$, the stochastic bounding inequality~\eqref{E:V_Noise_INEQ_LB} 
must not reduce to its deterministic counterpart~\eqref{E:LowerBound}. This 
requires that the diffusive term $\e \del \cdot(\vec{D} \del V_l)$ remains 
commensurate with the other terms in~\eqref{E:V_Noise_INEQ_LB}, at least near 
the fixed point at the origin.

The term $\e \del \cdot(\vec{D} \del V_l)$ can remain $O(1)$ near the origin as 
$\e\to0$ if $V_l$ develops a \emph{boundary layer} there. Purely polynomial 
$V_l$ can have no such boundary layer, so we add a non-polynomial term and let
\begin{equation}
\label{E:LogAnsatzV}
V_l(\x) = \alpha  \log[\e + \zeta(\x)] + P(\x).
\end{equation}
Here, $P(\x)$ is a polynomial, $\alpha$ is a tunable constant, and $\zeta(\x)$ 
is a positive definite quadratic form,
\begin{equation}
\zeta(\x) = \x^T \vec{Z} \x%,\quad \vec{Z}\succ \vec{0},
\end{equation} 
for a suitable symmetric, positive definite matrix $\vec{Z}$ (denoted 
$\vec{Z}\succ \vec{0}$) to be determined.  When $\e$ is small, the logarithmic 
term dominates $V_l$ near $\x=\vec{0}$, but $P(\x)$ remains significant away 
from the origin.

Despite not being polynomial, the ansatz~\eqref{E:LogAnsatzV} is suitable to 
derive a SoS optimization problem for the bound because the relevant derivatives 
of $V_l$ are rational functions,
\begin{subequations}
\begin{align}
\del V_l &= \frac{\alpha \del\zeta}{\e + \zeta} + \del P,	\\
\del\cdot(\vec{D}\del V_l) &= \alpha \frac{\del\cdot(\vec{D}\del\zeta)}{\e + 
\zeta} - \alpha\frac{ \del\zeta\cdot(\vec{D}\,\del\zeta)}{(\e + \zeta)^2} + 
\del\cdot(\vec{D}\del P).
\end{align}
\end{subequations}
Substituting these expressions into the bounding 
inequality~\eqref{E:V_Noise_INEQ_LB} gives a rational inequality with the 
positive denominator $(\e+\zeta)^2$. Multiplying this rational inequality by the 
denominator we obtain the equivalent polynomial inequality
\begin{equation}
\label{E:ineq_asymptotic}
\mathcal{L}(\x) := \mathcal{L}_0(\x) + \e\mathcal{L}_1(\x) + \e^2 
\mathcal{L}_2(\x) + \e^3\mathcal{L}_3(\x) \geq 0,
\end{equation}
where
\begin{subequations}
\begin{align}
\label{E:L0} 
%\notag
\mathcal{L}_0 &= \alpha \zeta \left( \f\cdot\del \zeta \right) + \zeta^2 \left( 
\f\cdot\del P + \phi - L \right) ,\\
\label{E:L1} 
%\notag
\mathcal{L}_1 &=	\alpha \zeta\, \del\cdot(\vec{D}\del\zeta) - \alpha 
\del\zeta\cdot(\vec{D}\,\del\zeta)  + \zeta^2\,\del\cdot(\vec{D}\del P) 
\\
&\quad + \alpha \f\cdot \del\zeta + 2\zeta\left(\f\cdot \del P +\phi - L\right)	
,\\
\label{E:L2} 
\mathcal{L}_2 &= \alpha \del\cdot(\vec{D}\del\zeta) + 2\zeta\, 
\del\cdot(\vec{D}\del P)	+ \f\cdot\del P + \phi - L	,\\
\label{E:L3} 
%\notag
\mathcal{L}_3 &= \del\cdot(\vec{D}\del P)	.	
\end{align}
\end{subequations}
Consequently, a lower bound $L$ on $\expect{\phi}{\e}$ can be calculated for a 
fixed, small $\e$ by solving the optimization problem
\begin{equation}
\label{E:OptLowerBound_Noise_NumericalLog}
\begin{aligned}
\max_{P,\alpha,\vec{Z}} \quad &L		\\
\text{s.t.}\quad  &\mathcal{L}(\x)  \in \Sigma, \\
%&\zeta(\x) = \x^T \vec{Z} \x, \quad 
&\vec{Z}\succ \vec{0}.
\end{aligned}
\end{equation}

The above SoS program can produce lower bounds larger than zero for very weak 
noise strengths, meaning that these bounds are not spoiled by the fixed point at 
the origin. In fact, while the SoS bounding program~\eqref{E:OptBoundFixedNoise} 
for polynomial $V_l$ reduces to its deterministic counterpart~\eqref{E:OptBound} 
as $\e\to0$, the above program does not. Instead, the polynomial $\mathcal L$ 
reduces to $\mathcal L_0$, retaining the term $\alpha \zeta \left( \f\cdot\del 
\zeta \right)$ that does not appear in the deterministic program. This term is 
due to the boundary layer built into $V_l$ around the origin.

\paragraph{Choosing the quadratic form $\zeta$.}

If the coefficients of the positive quadratic form $\zeta$ (equivalently, the 
entries of $\vec{Z}$) are tuned at the same time as the other coefficients, then 
the SoS constraint $\mathcal L\in\Sigma$ is quadratic in the decision variables, 
as opposed to linear. This produces a non-convex optimization problem for $L$ 
that is harder than a standard SDP and requires a bilinear solver. Here we keep 
the problem convex by fixing a non-optimal $\zeta$ in advance and tuning only 
the other coefficients in $V_l$. Even a non-optimal $\zeta$ provides a good 
lower bound at small $\e$ provided that, for reasons explained in 
\S\ref{SS:VanishingNoise}, we have not only $\zeta>0$ but also 
$\alpha\dot\zeta>0$ near the unstable fixed point. Quadratic $\zeta$ satisfying 
these two conditions 
%near an unstable fixed point 
can always be found if the point is a repeller but not if it is a saddle, which 
is why we have restricted ourselves to repellers. Further discussion of how to 
choose or optimize $\zeta$ is given in \S\ref{SS:VanishingNoise} and 
Appendix~\ref{S:ConstructZeta}.

%%%%%%%%%%%%%%%%%%%%
\paragraph{Example\,{\rm{:}} stochastic van der Pol oscillator with weak noise.}
\label{SS:ExampleNoisyVDPwithLog}

\begin{figure}[b!]
\centering
\subfigure[]{\label{F:LB_Noisy_vdp_Log} 
\includegraphics[width=0.47\textwidth, trim=0cm 3.5cm 5cm 0.5cm]
{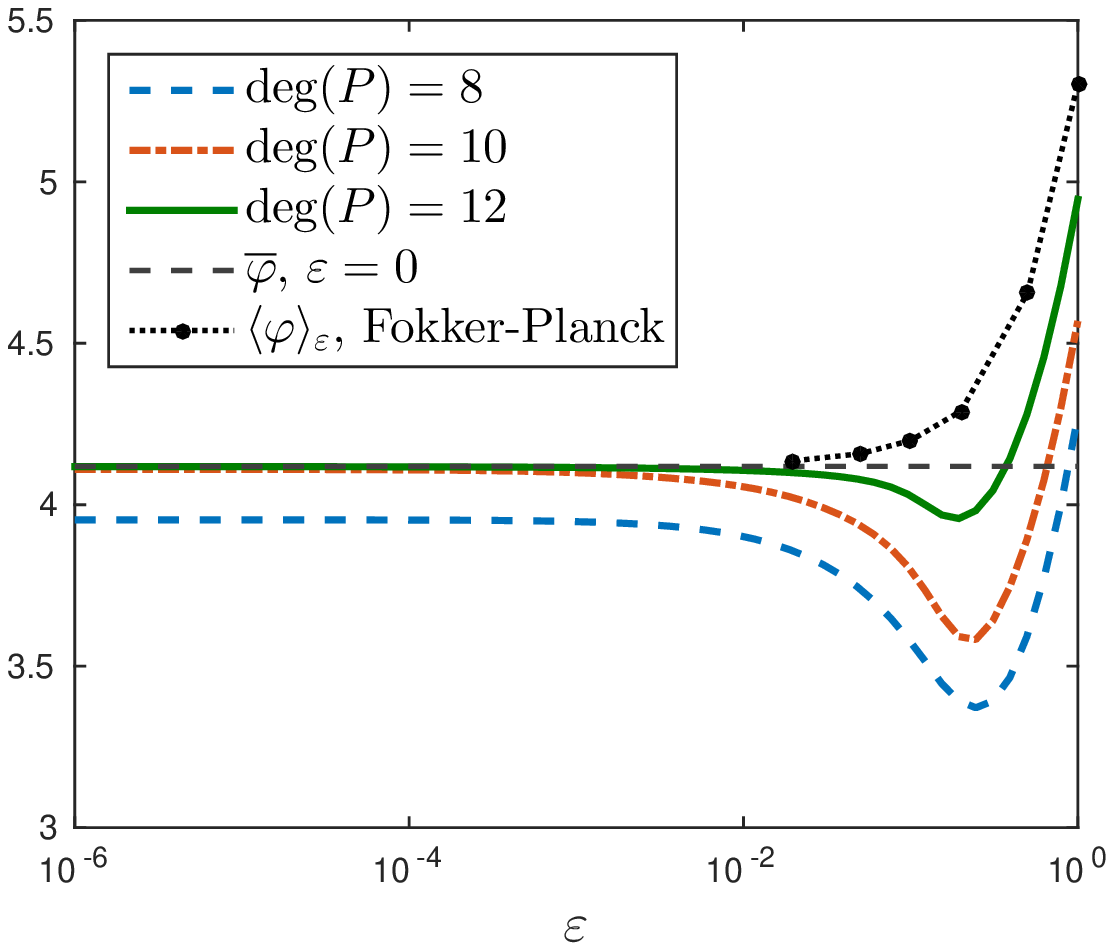} 
}
\subfigure[]{\label{F:LB_Noisy_LogvsNoLog} 
\includegraphics[width=0.47\textwidth, trim=0cm 3.5cm 5cm 0.5cm]
{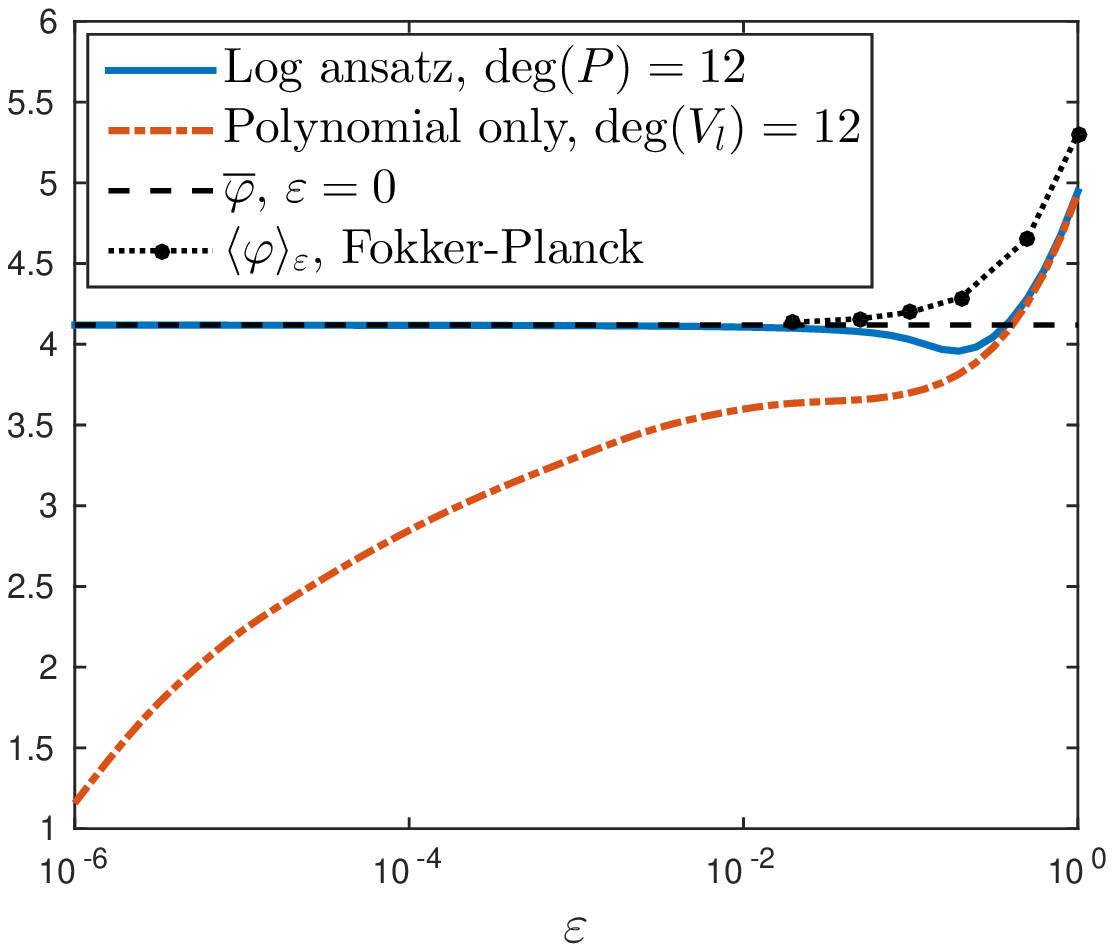}
}
\caption{{Lower bounds on $\expect{\phi}{\e} = \expect{x^2+y^2}{\e}$ for 
the stochastic van der Pol oscillator with $\mu=1$ 
and various noise strengths, obtained 
with~\eqref{E:OptLowerBound_Noise_NumericalLog} fixing $\zeta=x^2-xy+y^2$. (a) 
Bounds for different degrees of $P$, the polynomial part of $V_l$. (b) 
Comparison between the bounds computed with~\eqref{E:OptBoundFixedNoise} 
and~\eqref{E:OptLowerBound_Noise_NumericalLog} using polynomials of degree 
$12$. The deterministic ($\e=0$) 
average $\overline{\phi}=\overline{x^2+y^2}\approx 4.118$ 
is also shown for comparison.}}
\label{F:NoisyVDP_Log}
\end{figure}

We have used the methods of this section to derive lower bounds on 
$\expect{x^2+y^2}{\e}$ for the stochastic van der Pol 
oscillator~\eqref{E:VDP_Noise}. In the logarithmic ansatz~\eqref{E:LogAnsatzV} 
for $V_l$, we used the positive quadratic form $\zeta=x^2-xy+y^2$ that satisfies 
$\dot\zeta>0$ near the origin. Figure~\ref{F:NoisyVDP_Log} shows the resulting 
lower bounds computed using the SoS 
program~\eqref{E:OptLowerBound_Noise_NumericalLog} for  various degrees of the 
polynomial $P$. When the noise amplitude is moderate or small, these bounds 
dramatically improve on the bounds of Figure \ref{F:LB_Noisy_vdp} produced using 
purely polynomial $V_l$. When $\e\lesssim10^{-3}$, $L$ is indistinguishable from 
the deterministic \emph{local} lower bounds found in \S\ref{SS:LocalBounds}.

%%%%%%%%%%%%%%%%%%%%%%%% 
%  BOUNDS IN THE VANISHING NOISE LIMIT
%%%%%%%%%%%%%%%%%%%%%%%%
\section{Bounds on local attractors: the vanishing noise limit}
\label{SS:VanishingNoise}

In this section we propose a second method for bounding deterministic averages 
on a local attractor. While in \S\ref{SS:LocalBounds} we did this by segmenting 
the phase space, here we do it by adding noise. Many dynamical systems are 
stochastically stable, in the sense that the 
vanishing noise limit of a stationary expectation, $\lim_{\e\to0}\expect\phi\e$, 
is equal to the 
deterministic average $\overline\phi$ on a particular attractor~\cite{Young2002,Cowieson2005}. 
%Although we do not provide a formal proof, 
This is true in 
the van der Pol example, where the vanishing noise limit of $\expect{x^2+y^2}\e$ 
is equal to $\overline{x^2+y^2}$ on the limit cycle. Such correspondence can be 
exploited to bound $\overline\phi$ on a local attractor: by proving bounds on 
$\expect\phi\e$ that hold as $\e\to0$, we obtain bounds that hold also for 
$\overline\phi$. In essence, the vanishing noise limit preserves the attractor 
while destroying unstable invariant structures. This idea was proposed 
in~\cite{Chernyshenko2014a}, though here we extend it by treating the $\e\to0$ 
limit rigorously, rather than numerically.

We continue the analysis of~\S\ref{SS:SmallNoise}, still assuming that $\x=\vec{0}$ is a repelling fixed point, $\phi(\vec{0})=0$, 
and $\overline{\phi}>0$ on the attractor of interest. We suppose also that 
$\lim_{\e\to0}\expect\phi\e = \overline\phi$, meaning that the zero-noise limit of the expectation converges to the deterministic time average. Rigorous statements about when this property holds can be found in~\cite{Young2002,Cowieson2005} and references therein.

Under these assumptions, we can obtain a lower bound $L\le\overline\phi$ 
that is tight for the local attractor by deriving a lower 
bound $L\le\lim_{\e\to0}\expect{\phi}{\e}$.
To this end, we recall that although the polynomial inequality $\mathcal L\ge0$ 
considered in~\S\ref{SS:SmallNoise} suffices to prove the lower bound, we are 
really interested in proving the integral inequality~\eqref{E:LB_IntIneq}. 
It is proven in Appendix~\ref{S:AnalyticEpsProof}, under mild assumptions on the 
stationary distribution $\rho$, that~\eqref{E:LB_IntIneq} holds in the limit 
$\e\to0$ if there exists any $\gamma>0$ such that $\mathcal{L}_0\geq\gamma 
\zeta^2$. Dividing this relation by $\zeta$ (which is constructed to be positive definite) and 
rearranging gives the sufficient condition
\begin{equation}
\label{E:VanishingNoiseSuffCond}
\alpha\, \f\cdot\del \zeta + \zeta \left[ \f\cdot\del P + \phi - (L-\gamma) 
\right]  \geq 0.
\end{equation}
Since $\gamma>0$ decreases $L$ by an arbitrarily small amount, we can let the 
bound be implied by a strict inequality and set $\gamma=0$. The lower bound we 
seek, $L<\lim_{\e\to0}\expect{\phi}{\e}$, is thus returned by the SoS 
optimization problem
\begin{equation}
\label{E:SoSproblem}
\begin{gathered}
\max_{P,\alpha,\vec{Z}} \quad L
\\[5pt]
\begin{gathered}
\text{s.t.} \quad
\alpha\, \f\cdot\del \zeta  + \zeta \left( \f\cdot\del P + \phi - L \right) \in 
\Sigma,\\
%\zeta(\x) = \x^T \vec{Z} \x, \quad 
\vec{Z}\succ \vec{0}.
\end{gathered}
\end{gathered}
\end{equation}

Like the finite-noise program~\eqref{E:OptLowerBound_Noise_NumericalLog}, the 
above program can be made convex by fixing in advance the positive definite 
quadratic form $\zeta$. The chosen $\zeta$ must satisfy $\alpha\dot\zeta>0$ near 
the origin. That is, it must satisfy 
\begin{equation}
\label{E:zetaCond}
\alpha\,\tilde{\f}\cdot\del\zeta > 0,
\end{equation}
where $\tilde{\f} = \vec{J}_0\,\x$ denotes the linearized dynamics near the 
origin. This condition is needed because the SoS constraint becomes
\begin{equation}
\label{E:ReducedBalanceCondition}
\alpha \, \tilde{\f}\cdot\del \zeta - L \zeta \gtrsim 0
\end{equation}
near the origin. The only way the above condition can be satisfied for $L$ 
larger than zero is if its first term is positive. 
Appendix~\ref{S:ConstructZeta} gives a way of constructing an admissible $\zeta$ 
when the unstable fixed point is repelling. When the unstable fixed point is a 
saddle, it is not generally possible to satisfy~\eqref{E:zetaCond} for positive 
definite $\zeta$.

%%%%%%%%
\paragraph{Example\,{\rm{:}} deterministic bounds for the van der Pol limit 
cycle---vanishing noise formulation.}

To demonstrate the methods of this section, we have computed lower bounds on the 
vanishing noise limit of $\expect{x^2+y^2}\e$ for the van der Pol system using 
the SoS program~\eqref{E:SoSproblem}. The results serve also as deterministic 
lower bounds on $\overline\phi$ for all trajectories approaching the limit 
cycle (strictly speaking, this conclusion requires proving that $\lim_{\e\to0}\langle x^2 + y^2 \rangle_\epsilon=\overline{x^2 + y^2}$, which is outside our present scope). Bounds of the latter type appear also in \S\ref{SS:LocalBounds}, computed 
differently using the SoS formulation~\eqref{E:OptBoundSproc}.
\begin{table}%[b!]
\footnotesize	% required by SIAM!
\caption{Choices of $\zeta$ for the van der Pol oscillator. }
\centering
\begin{tabular}{c | c  c}
\hline \hline
 & $\mu\leq 2$ & $\mu > 2$ \\
 \hline\\[-8pt]
 $\zeta_1$ & $x^2- xy + y^2$ & $x^2-xy + y^2$ \\
 $\zeta_2$ & $x^2-\mu xy + y^2$ & $\mu x^2-4xy + \mu y^2$  \\
 $\zeta_3$ & $x^2-\mu xy + y^2$ & $(\mu^2-2)x^2-2\mu xy + 2y^2$\\ \hline
\end{tabular}
\label{Table:Table1}
\end{table}
\begin{figure}%[b!]
\centering
\label{F:LB_ZeroNoise_vdp} 
\includegraphics[width=0.57\textwidth, trim=0cm 0.75cm 0cm 0cm]
{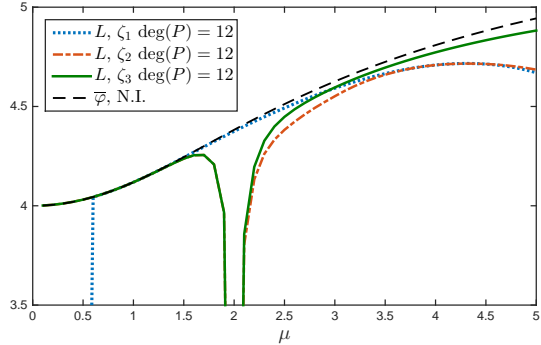}
\caption{{Lower bounds on $\overline{\phi}=\overline{x^2+y^2}$ for 
the van der Pol oscillator obtained 
with the vanishing noise program~\eqref{E:SoSproblem} 
for $P$ of degree 12, along with the values of $\overline{\phi}$ 
found by numerical integration (N.I.).}}
\label{F:ZeroNoiseVDP}
\end{figure}
Figure~\ref{F:ZeroNoiseVDP} shows the bounds obtained by fixing the degree of 
the polynomial $P$ to 12 for the three different $\zeta$ defined in 
Table~\ref{Table:Table1}. The condition~\eqref{E:zetaCond} is satisfied with 
$\alpha>0$ by $\zeta_2$ and $\zeta_3$ when $\mu\neq2$, and by $\zeta_1$ when 
$4-2\sqrt{3}<\mu<4+2\sqrt{3}$. The failure of the condition is why the bound 
using $\zeta_1$ is poor for $\mu\leq 4-2\sqrt{3} \approx 0.54$, while the bounds 
using $\zeta_2$ or $\zeta_3$ are poor near $\mu=2$. Taking the best of the three 
lower bounds at each $\mu$ gives a lower bound on $\overline{x^2+y^2}$ within 
1\% of the value found by numerical integration over the entire range $0.1\le\mu\le5$.

%%%%%%%%%%%%%%%%%%%%%%%% 
%  FUTURE DIRECTIONS
%%%%%%%%%%%%%%%%%%%%%%%%
\section{Future directions}
\label{S:FurtherComments}

Despite our success in computing bounds for the van der Pol oscillator, the 
study of more complicated dynamical systems presents several obstacles that 
require further investigations. We discuss some of these below, and put forward 
some preliminary ideas for future improvements.

\subsection{Bounds for high-dimensional systems}
\label{SS:HighDimBounds}
The techniques we have developed in this paper can in principle be applied to 
systems of arbitrarily high but finite dimension. However, computing tight bounds 
typically requires the use of polynomial storage functions of large degree.
Since the size of the SDP required to determine the 
existence of a SoS decomposition for a polynomial of degree $d$ in $n$ variables 
grows as $\left(\begin{smallmatrix}n+d\\d\end{smallmatrix}\right)$~\cite[Theorem 
3.3]{Parrilo2003}, relatively large SDPs must be solved even for a 
low-dimensional system such as the van der Pol oscillator. Moreover, the SDPs 
seem to be consistently ill-conditioned, and more so for larger systems. Our present methods are therefore practical only for systems of moderate dimension.

%%%% ORIGINAL
%One way to improve the numerical conditioning of the SDPs is to try to rescale 
%the system, although a systematic rule does not exist. A preliminary 
%investigation showed that rescaling the van der Pol system such that the limit 
%cycle is contained in the box $[-1,1]^2$ dramatically improves the numerical 
%conditioning, and bounds indistinguishable from those presented in the previous 
%sections can be obtained efficiently using standard solvers.

%%% DAVID - IMPROVED!
One way to improve the numerical conditioning of the SDPs is to try to rescale 
the system. A preliminary investigation showed that rescaling the van der Pol 
system such that the limit cycle is contained in the box $[-1,1]^2$ dramatically 
improves the numerical conditioning, and bounds very similar to those reported 
above can be obtained without the need for multiple precision solvers. 
Researchers studying other aspects of dynamics using SDPs have had similar 
success rescaling the relevant dynamics to lie in $[-1,1]^n$. However, the 
appropriate rescaling is not generally clear \emph{a priori}.

%%%% ORIGINAL
%The computational cost of solving large SDPs, instead, could be reduced by 
%exploiting special structures in the optimization problems, such as symmetries 
%and/or any sparsity in the data, for instance using the ideas 
%of~\cite{Fukuda2000,Nakata2003}. Although currently unknown, it may also be 
%possible that special structures in the system dynamics (e.g. Hamiltonian) can 
%be exploited to assist the formulation of more specific ansatze for the storage 
%functions, thereby reducing the computational burden of the optimization. 
%Finally, an alternative approach might be to work with a 
%system of smaller dimension, obtained truncating a high 
%dimensional system and bounding the errors introduced by the 
%truncation~\cite{Goulart2012,Chernyshenko2014a,Huang2015}.

%%% DAVID - IMPROVED!
There are several possible ways to reduce the computational cost of constructing bounds. One option is to exploit special structures in the SDPs such as symmetries or sparsity (e.g.~\cite{Fukuda2000,Nakata2003}). It might also be profitable to exploit special structure of the underlying dynamics, such as conserved quantities. Finally, it is sometimes possible to work with a lower-dimensional truncation of a high-dimensional system and bound the errors introduced by the truncation (e.g.~\cite{Goulart2012,Chernyshenko2014a,Huang2015}).

\subsection{Bounds for systems governed by partial differential equations}

% DAVID - SINGLE PARAGRAPH
The methods developed here apply only to finite-dimensional systems, and more 
work is needed to extend them to partial differential equations (PDEs). One idea 
is to project the PDE variables onto a finite Galerkin basis. Bounds can then be 
constructed using the finite-dimensional system, provided that the influence of 
the unprojected component can be controlled. This approach has been used 
successfully for nonlinear stability analysis of a fluid dynamical 
PDE~\cite{Goulart2012,Chernyshenko2014a,Huang2015}. A second idea is to bound 
integrals of the PDE variables directly using the dissipation inequalities 
proposed in~\cite{Ahmadi2016}. In essence, these inequalities are integral 
constraints that replace the polynomial constraints of 
Proposition~\ref{T:Proposition1} when a finite-dimensional system is replaced by 
a PDE. Despite recent advances in the numerical implementation of integral 
inequalities~\cite{Valmorbida2016}, however, this technique has so far been 
applied only to simple examples.

% ORIGINAL
%Since many physical systems are governed by partial differential equations 
%(PDEs), it is in the interest of future work to derive techniques to compute 
%bounds on solution to PDEs.  In light of \S\ref{SS:HighDimBounds}, the 
%infinite-dimensional nature of PDEs poses a challenge for a direct extension of 
%our techniques, even in the simplest deterministic case. 
%
%As a first approach, deterministic PDEs could be studied with the techniques 
%developed in this work via a rigorous projection onto a finite-dimensional space 
%of low to moderate dimension. Similar ideas have already been applied 
%in~\cite{Goulart2012,Chernyshenko2014a,Huang2015} to study the stability of 
%fluid flows. However, such rigorous projections typically rely on (conservative) 
%estimates for the unknown solutions of a PDE, and therefore the optimal bounds 
%computed using Propositions~\ref{T:Proposition1} and~\ref{T:Proposition2} are 
%not likely to be tight. 
%
%Alternatively, deterministic PDEs could in principle be analyzed directly using 
%the \emph{dissipation inequalities} proposed in~\cite{Ahmadi2016}. In essence, 
%these are integral inequality constraints that, in some sense, replace the 
%inequality constraints of Proposition~\ref{T:Proposition1} when a 
%finite-dimensional system is replaced by a PDE. Despite recent advances in the 
%numerical implementation of integral inequalities~\cite{Valmorbida2016}, 
%however, this technique has only been applied to simple examples, and the 
%application to problems of practical interest remains a challenge for future 
%work.

\subsection{Tight bounds for systems with saddle points}

If deterministic bounds are spoiled by a saddle point $\x_s$ that is not 
embedded in the attractor of interest, the \Sprocedure~can be used to remove a 
set containing $\x_s$ that is disjoint from the attractor. 
However, it can be difficult to establish that $\x_s$ is indeed separate
from the attractor. 
If $\x_s$ \emph{is} embedded in the attractor, alternative methods 
must be found. 
The vanishing-noise approach of~\S\S\ref{SS:SmallNoise}--\ref{SS:VanishingNoise} seems promising, but the present formulation works only for repelling points, not saddle points.
This is because~\eqref{E:ReducedBalanceCondition} requires that the polynomial 
$\alpha\zeta(\x)$ increases along all trajectories near the unstable fixed point, but 
$\alpha\zeta$ is monotonic and cannot increase along trajectories on the stable 
and unstable manifolds of $\x_s$ simultaneously.  

Similarly, the logarithmic ansatz~\eqref{E:LogAnsatzV} is not suitable when 
stochastic bounds with weak (but finite) noise are spoiled by a saddle point. 
In fact, the term $\f\cdot\del V$ in each of the inequalities 
in~\eqref{E:StochasticInequalities} has opposite signs along the stable and 
unstable manifolds of $\x_s$ unless $\del V$ changes rapidly, and any negative 
contribution must be balanced by large Laplacians. This requires polynomials of 
large degree, making the SoS programs intractable in practice. 

One possible solution is to find an ansazt for $V$ that, similarly 
to~\eqref{E:LogAnsatzV} for repelling fixed points, lets stochastic bounds stay 
tight as $\e\to0$ by developing appropriate boundary layers.
The following proposition suggests a way to deduce the ideal scaling of $V$ in 
this limit.
\begin{proposition}
\label{T:TheoremSC}
Consider the inequality
\begin{equation}
S(\x) = \e \del \cdot(\vec{D} \del V_l) + \f\cdot\del V_l + \phi - L \geq 
0\qquad  \forall \x \in \Real^n \label{E:Sineq}
\end{equation}
that is a sufficient condition for the lower bound $L\le\expect{\phi}{\e}$. If 
$S(\x)$ is continuous and the stationary distribution $\rho(\x)$ is piecewise 
continuous, then the above inequality can be satisfied for the perfect lower 
bound $L=\expect{\phi}{\e}$ only if $S(\x)\equiv0$ wherever $\rho(\x)>0$.
\end{proposition}

This statement is proven by adding $\expect{\phi-L}\e=0$, which holds because 
$L$ is a perfect bound, and expression~\eqref{E:zero} to find $\expect S\e=0$. 
Since $S(\x)$ is continuous and non-negative, $\expect S\e=0$ is possible only 
if $S(\x)=0$ wherever $\rho(\x)>0$, as claimed. An analogous statement holds for 
the perfect upper bound.

Proposition~\ref{T:TheoremSC} suggests that as $V_l$ is improved and $L$ gets 
closer to its perfect value of $\expect\phi\e$, the inequality $S(\x)\ge0$ gets 
closer in some sense to being an equality. This motivates us to consider 
solutions $V(\x)$ to the equation
\begin{equation}
\label{E:AsymptitocEqn}
\e \del \cdot(\vec{D} \del V) + \f\cdot\del V + \phi - L=0,
\end{equation}
where $L=\expect\phi\e$ and $\rho$ obeys the decay 
condition~\eqref{E:boundaryInt}. A solution exists only when $L = 
\expect{\phi}{\e}$, as seen by taking the expectation of the equality and 
using~\eqref{E:zero}. Asymptotic analysis of $V(\x)$ in the above equality as 
$\e\to0$ might suggests a good ansatz for $V_l$ in the corresponding inequality.

\subsection{Assessing the quality of numerical bounds}

Proposition~\ref{T:TheoremSC} also provides a good way of assessing \emph{a 
posteriori} whether the stochastic bounds obtained with SoS optimisation are 
nearly sharp. This is useful, for instance, when data from numerical 
simulations are not available for 
comparison. For example, if a lower bound is close to the exact value of 
$\expect\phi\e$, $V_l$ should be close to the solution 
of~\eqref{E:AsymptitocEqn} whenever $\rho(\x)\neq0$. In addition, an argument 
similar to the proof of Proposition~\ref{T:TheoremSC} for deterministic bounds 
shows that if an upper or lower bound is perfect, then
\begin{equation}
\label{E:CorrectBehaviourCheck}
\dot{V}=\overline{\phi}-\phi(x)
\end{equation}
for any trajectory $\x(t)$ on the attractor. If a good lower bound is proven 
using $V_l$, for instance, then $\dot V_l$ should be very close to 
$\overline{\phi}-\phi(x)$ on the attractor. 

\begin{figure}%[t!]
\centering
\includegraphics[scale=0.6, trim = 0cm 0.5cm 0cm 0.5cm]
{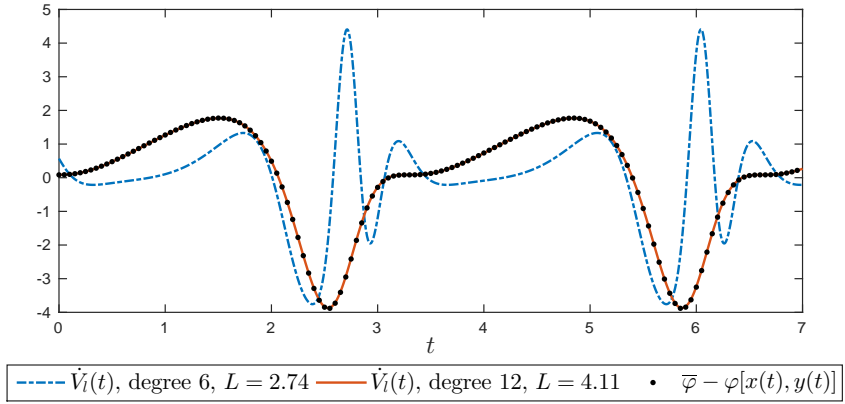}
\caption{Comparison of $\dot{V_l}$ along the limit cycle of the deterministic 
van der Pol oscillator for $\mu=1$ 
(where $\overline{\phi}=\overline{x^2+y^2}\approx 4.118$). The 
logarithmic ansatz~\eqref{E:LogAnsatzV} for $V_l$ was used with $P$ of order 6 
and 12. The storage functions were calculated with~\eqref{E:SoSproblem}, while 
the limit cycle is from numerical integration.}
\label{F:CorrectBehaviour} 
\end{figure}

For illustration, we have solved the lower bound program~\eqref{E:SoSproblem} 
for the van der Pol oscillator with $P(\x)$ of degree 6 and degree 12. The 
degree-6 lower bound of $L=2.74$ is rather poor, while the degree-12 bound of 
$L=4.11$ is nearly tight. Figure~\ref{F:CorrectBehaviour} shows the variation of 
$\dot{V_l}(t)$ along the limit cycle in each case and compares it to the value of 
$\overline{\phi}-\phi(x)$ along the limit cycle. The $\dot{V_l}(t)$ and 
$\overline{\phi}-\phi(x)$ curves match very closely in the degree-12 case but 
not in the degree-6 case.

\subsection{Formulating computer-assisted proofs}
\label{SS:ComputerAssistedProof}

The bounds we have computed for the van der Pol system all rely on the numerical 
solution of SDPs, so they are subject to roundoff error.  It is unlikely that 
roundoff error has led to any invalid bounds in our calculations because we have 
solved the SDPs using multiple-precision arithmetic and checked the results by 
increasing the precision. 
%% ORIGINAL
%Nonetheless, such computations cannot be considered 
%rigorous to the standard of computer-assisted proofs unless bounds are proven on 
%the effects of roundoff errors. This requires \emph{a posteriori} verification 
%of the numerical solution returned by an SDP solver. Theorem 4 
%in~\cite{Lofberg2009} provides a sufficient condition for bounding numerical 
%errors in SoS decompositions. Appendix~\ref{S:RigorousBounds} demonstrates the 
%use of this condition and offers further discussion.
%% DAVID - IMPROVED
Nonetheless, such computations cannot be considered rigorous to the standard of computer-assisted proofs because of the presence of roundoff errors. Two different methods have been proposed for obtaining rigorous results from numerical SDP solutions. The first method is to project an approximate numerical solution onto an exact solution in terms of rational numbers~\cite{Peyrl2008,Kaltofen2012}. The second method is to use perturbation analysis, made rigorous by interval arithmetic, to construct a small interval around the approximate numerical solution in which the exact solution is guaranteed to lie. This approach has been implemented in the software VSDP~\cite{Jansson06}, but the relevant SDPs would need to be constructed manually since there is no sum-of-squares parser available that incorporates interval arithmetic.

%%%%%%%%%%%%%%%%%%%%%%%% 
%  CONCLUSIONS
%%%%%%%%%%%%%%%%%%%%%%%%
\section{Conclusions}
\label{S:Conclusion}

In this work we have presented computer-assisted methods for deriving bounds on 
average quantities in both deterministic and stochastic dynamical systems using 
sum-of-squares programming. We have given particular attention to proving bounds 
that apply only to trajectories approaching a particular attractor. 
One method is to use the \Sprocedure~to omit segments of phase space on which 
the bounds do not need to hold. Another strategy is to remove unstable invariant 
structures by adding noise to the system. This idea, proposed previously 
in~\cite{Chernyshenko2014a}, has been extended here by analyzing the weak and 
vanishing noise cases when the unstable structures to be omitted are repelling 
fixed points. These methods give improved bounds for weak but 
finite noise, and also in the rigorous limit of vanishingly weak noise.

\begin{figure}[t!]
\centering
\subfigure[]{
\includegraphics[width=0.47\textwidth, trim=0cm 0cm 0.95cm 0cm]
{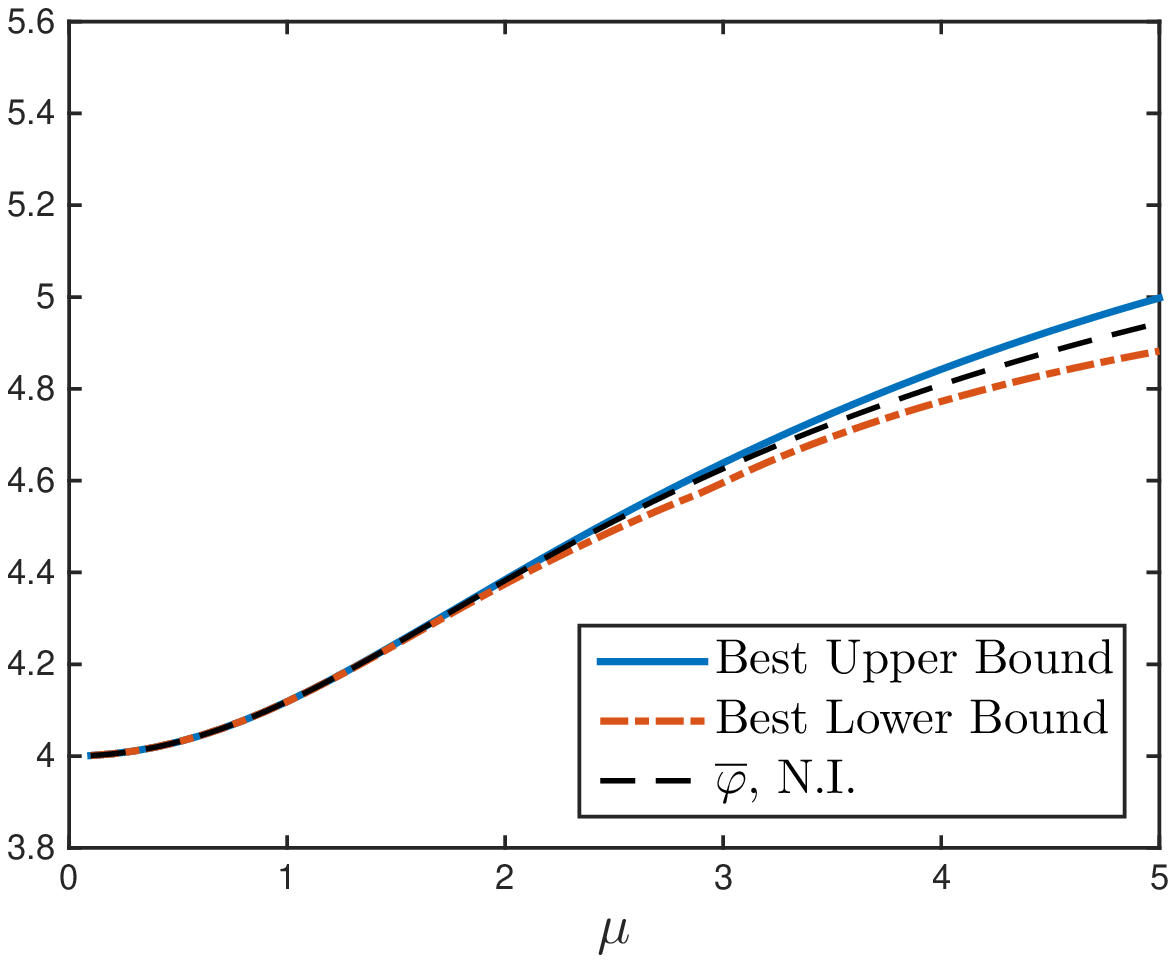}
}
\subfigure[]{
\includegraphics[width=0.47\textwidth, trim=0.95cm 0cm 0cm 0cm]
{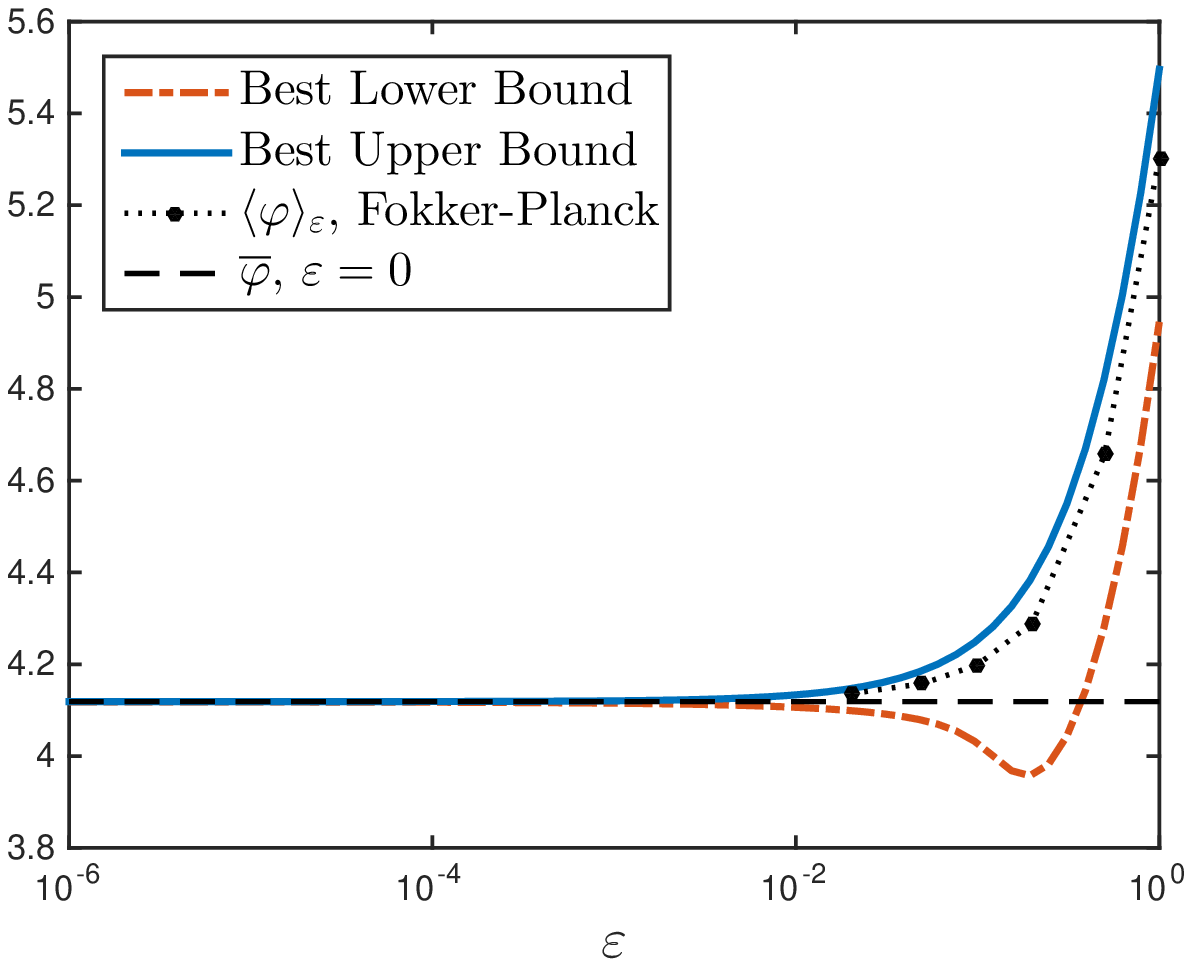}
}
\caption{Best bounds obtained on averages of $\phi=x^2+y^2$ for the van 
der Pol oscillator. (a) Bounds on the 
deterministic time average $\overline{\phi}$ over the limit cycle for various 
damping strengths $\mu$. (b) Bounds on the stochastic expectation 
$\expect{\phi}{\e}$ with $\mu=1$ and various noise strengths $\e$.}
\label{F:BestBounds_vdp}
\end{figure}

Our methods have worked well when applied to the van der Pol oscillator. The 
best deterministic and stochastic bounds proven throughout are summarized in 
Figure \ref{F:BestBounds_vdp}. In the deterministic case, we obtained upper and 
lower bounds on the infinite-time average of ${x^2+y^2}$ over the limit cycle 
that are all within 1\% of the ``true'' values found by numerical integration. 
In the stochastic case, bounding the stationary expectation of $x^2+y^2$ at a 
variety of noise strengths, we obtained upper bounds all within 1\% and lower 
bounds all within 10\% of the ``true'' values found by solving the Fokker-Planck 
equation numerically.

Moving to dynamical systems other than the van der Pol oscillator, whether the 
methods we have described can yield similarly tight deterministic and stochastic 
bounds depends on the details of those systems. If the dimension of a system is 
small enough for SoS optimization to be computationally feasible,
we expect that tight global bounds on deterministic averages can be obtained 
using the methods of \S\ref{SS:SoSBoundsReview}. In fact, since the first draft of this work, the bounding techniques we have presented for deterministic systems have been successfully applied to the design of control systems for fluid flows~\cite{Huang2016,Lasagna2016}.
With stochastic forcing that is 
not too weak, we expect that fairly tight bounds on stationary expectations also 
can be obtained using the methods of \S\ref{SS:GeneralNoiseFormulation}. As the 
noise strength decreases, the tendency of unstable invariant structures to spoil 
tight bounds can be combatted by the methods of \S\ref{SS:SmallNoise} if these 
structures are repelling fixed points. New methods will be needed to maintain 
tight bounds as noise becomes weak in systems with other unstable structures, 
including saddle points. Lastly, the methods of \S\ref{SS:LocalBounds} can be 
used to bound averages on trajectories within an absorbing domain and, in 
particular, on a single attractor. This much is true for any dynamical system, 
but whether these bounds can be tight depends on the nature of the attractor. 
For instance, bounds for a chaotic attractor must apply not just to generic 
chaotic trajectories but also to all trajectories embedded within the attractor, 
such as saddle points and unstable orbits. What can be deduced about chaotic 
attractors is of particular interest if the methods presented here are to be 
applied to complex systems of physical and engineering relevance.

%%%%%%%%%%%%%%%%%%%%%%%% 
%    APPENDIX
%%%%%%%%%%%%%%%%%%%%%%%%
\appendix

%%%%%%%%%%%%%%%%%%%%%%%%
\section{Towards rigorous bounds}
\label{S:RigorousBounds}
As explained in \S\ref{SS:ComputerAssistedProof}, techniques are available to control the roundoff errors in the solution of the SDPs and compute rigorous bounds. However, these may be impractical for large problems. We are also unaware of any parsers for SoS programs that incorporate rigorous computations.  For these reasons, we illustrate a method to produce  rigorous bounds that can be implemented with existing software.  For definiteness, we consider the upper bound problem in~\eqref{E:OptBound}. 

The solution of the SoS program consists of a polynomial storage function $V_u$, 
a bound $U$, a vector of monomials $\vec{z}(\x)$ and a positive semidefinite, 
symmetric matrix $\vec{Q}$ such that
\begin{equation}
\vec{z}(\x)^T\vec{Q}\vec{z}(\x) = U-\phi(\x)-\f\cdot\del V_u .
\end{equation}
Since $\vec{Q}$ is positive semidefinite, $\vec{z}(\x)^T\vec{Q}\vec{z}(\x)$ is a 
SoS polynomial and so is the right-hand side. We refer the reader 
to~\cite{Parrilo2003,Lofberg2009} for further details.

In practice, the decomposition is only approximate, and the error
\begin{equation}
e(\x) := \vec{z}(\x)^T\vec{Q}\vec{z}(\x) -\left[ U-\phi(\x)-\f\cdot\del V_u 
\right]
\end{equation}
is non-zero. However, according to Theorem 4 in~\cite{Lofberg2009}, the 
polynomial $U-\phi(\x)-\f\cdot\del V_u$ is still certifiably non-negative if
\begin{equation}
\label{E:RigorousChecks}
\lambda_0 - \text{dim}(\vec{Q})\times\abs{r} \geq 0,
\end{equation}
where $r$ is the coefficient of $e(\x)$ of largest magnitude, $\lambda_0$ is the 
smallest eigenvalue of $\vec{Q}$ and $\text{dim}(\vec{Q})$ denotes the dimension 
of $\vec{Q}$. 

Unfortunately, the optimal solution of a SoS problem rarely satisfies this 
condition, since $\lambda_0\approx0$ for the optimal 
$\vec{Q}$~\cite{Lofberg2009}. 
Consequently, we suggest to compute a slightly suboptimal $U$ by manually 
decreasing an initial guess $U_0$ using a sequence of feasibility problems, 
checking after each step that~\eqref{E:RigorousChecks} holds using rigorous 
computations such as interval arithmetic~\cite{Alefeld2000,Hladik2010}. The 
steps are outlined in Algorithm~\ref{Algorithm1}.
\begin{algorithm}[b!]
\caption{. Sequence of feasibility checks to compute rigorous SoS upper bounds}
\label{Algorithm1}
\small
\vskip 0.5ex
\begin{algorithmic}[1]
\State $\textit{U} \gets$ initial guess for the upper bound
\BState \emph{loop:}
\While {rigorous checks not verified \textbf{and} $U$ has not converged}
\State $\textit{V} \gets$  construct a polynomial with variable coefficients
\State $\textit{expr} \gets$ construct the expression $U - \f\cdot\del V_u - 
\phi$
\If {find a suitable $\textit{V}$ and a SoS decomposition 
$\textit{expr}=\vec{z}^T\vec{Q}\vec{z}$} 
	\State $V,\vec{Q} \gets$ round the coefficients of $V$ and the entries of 
$\vec{Q}$ to $d$ decimal places
	\State Check~\eqref{E:RigorousChecks} for the truncated $V$, $\vec{Q}$ using 
rigorous numerics 
	\If {checks are verified}
		\State decrease $U$ and \textbf{goto} \emph{loop}
	\Else
		\State increase $U$ and \textbf{goto} \emph{loop}
	\EndIf
\Else
	\State increase $U$ and \textbf{goto} \emph{loop}
\EndIf
\EndWhile
\end{algorithmic}
\end{algorithm}

\begin{figure}%[b!]
\centering
\includegraphics[width=0.6\textwidth, trim=0cm 0cm 0cm 0.5cm]
{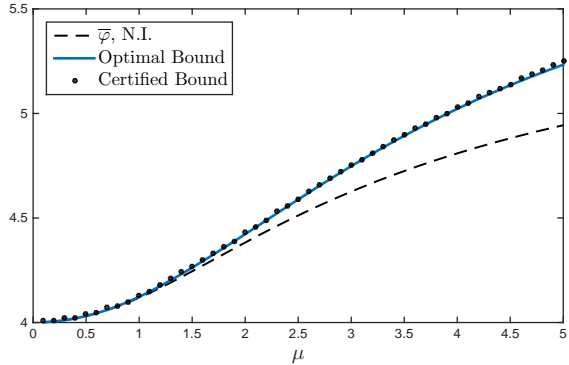}
\caption{{Comparison between the optimal deterministic upper bounds on 
$\overline{\phi}=\overline{x^2+y^2}$ of \S\ref{SS:ExampleDeterministicVDP_UB} 
and the certified bounds 
obtained with Algorithm~\ref{Algorithm1} for $V_u$ of degree 10. The ``true'' 
value of $\overline{\phi}$ from numerical integration (N.I.) is also shown.}}
\label{F:UBvdpRigorous}
\end{figure}
Figure~\ref{F:UBvdpRigorous} compares some preliminary bounds obtained with 
Algorithm~\ref{Algorithm1} for the van der Pol oscillator to the optimal ones.  
The results are only slightly suboptimal; moreover the feasibility problems were 
better conditioned than the full SoS optimization and could be solved with 
standard double-precision solvers such as SeDuMi~\cite{Sturm1999}. We remark 
that given the preliminary stage of this investigation, we avoided the technical 
step of computing rigorous roundoff errors for the eigenvalues of $\vec{Q}$. 
Instead, we checked condition~\eqref{E:RigorousChecks} using arbitrary precision 
arithmetic, increasing the precision until $\lambda_0$ changed by less than 1\%. 
Although our results are only ``certified'' but not fully rigorous, we expect 
that similar results would be obtained when $\lambda_0$ is computed rigorously.
Moreover, we also stress that Algorithm~\ref{Algorithm1} was not very robust and 
required careful tuning (e.g. in the choice of decimal figures to retain). A 
thorough investigation and resolution of these numerical issues is left for 
future work.
%

%%%%%%%%%%%%%%%%%%%%%%%%
\section{Construction of an absorbing domain for the van der Pol oscillator}
\label{S:VDPAbsorbingDomainProof}
To show that $\mathcal{T}_r = \set{(x,y) \, |\, x^2+y^2-r^2 \geq 0 }$ is an 
absorbing domain for the van der Pol oscillator for $r\leq1$, let us reverse the 
direction of time in~\eqref{E:VDP_SS} and consider the energy $E = x^2+y^2$ of 
the system $\dot{\x}=-\f(\x)$, for which the origin is a stable fixed point. One 
has
\begin{equation}
\begin{aligned}
\dot{E}= 2x\dot{x} + 2y\dot{y} 
		    = -2xy - 2\mu y^2(1-x^2) + 2xy 
		    = 2\mu y^2(x^2-1)
\end{aligned}
\end{equation}
meaning that $E\leq 0$ when $\abs{x}\leq1$. Therefore, any contour of $E$ which 
is contained in the strip $\abs{x}\leq1$ defines the boundary of a region of 
attraction of the origin for the time-reversed oscillator.  One concludes that 
in the original system all orbits (except of course the fixed point at the 
origin) will leave the ball $x^2+y^2 < r^2$ if $r\leq 1$, i.e. $\mathcal{T}_r$ 
is an absorbing domain.

%%%%%%%%%%%%%%%%%%%%%%%%
\section{Bounds with vanishing noise}
\label{S:AnalyticEpsProof}
To prove a lower bound on $\expect{\phi}{\e}$ we need to prove 
inequality~\eqref{E:LB_IntIneq}. The ansatz~\eqref{E:LogAnsatzV} for $V$ allows 
us to rewrite this condition as
\begin{equation}
\label{E:VanishingNoise_Thesis0}
\left\langle \frac{\mathcal{L}(\x)}{\left[ \e + \zeta(\x)\right]^2} 
\right\rangle_{\e} \geq 0,
\end{equation} 
where $\mathcal{L}$ is as in~\eqref{E:ineq_asymptotic}. Since $\zeta$ is 
quadratic in $\x$, it can be verified that $\mathcal{L}_0$ is the dominant term 
in $\mathcal{L}$ for any fixed $\x\neq 0$ in the limit $\e\to 0$. Moreover, 
$\mathcal{L}_0$ is the dominant term as $\vert\x\vert\to\infty$ because it 
contains the monomials of highest order. The condition 
$\mathcal{L}_0\geq\gamma\zeta^2$ then implies that $\mathcal{L}(\x)\geq 0$ in 
the limit $\e\to 0$ for all $\x$, with the possible exception of a  ball $B_R$ 
of radius $R\sim\e^{\eta}$, $0<\eta<1/2$, where the integrand 
in~\eqref{E:VanishingNoise_Thesis0} becomes singular.  For definiteness, we will 
write $R=C\e^{\eta}$ for some constant $C$.  Thus, we conclude that
\begin{equation}
\lim_{\e\to0}\int\limits_{\Real^n\smallsetminus 
B_R}\rho(\x)\frac{\mathcal{L}(\x)}{\left[ \e + \zeta(\x)\right]^2}d\x \geq 0.
\end{equation}
To prove that~\eqref{E:VanishingNoise_Thesis0} holds, we will now show that
\begin{equation}
\label{E:VanishingNoise_Thesis}
\lim_{\e\to0} \int\limits_{B_R}  \rho(\x)\frac{\mathcal{L}(\x)}{\left[ \e + 
\zeta(\x)\right]^2} d\x = 0.
\end{equation}
Upon substitution of~\eqref{E:LogAnsatzV}, the integral becomes
\begin{multline}
\int\limits_{B_R}   \rho \left[  
\e\,\frac{\alpha\del\cdot(\vec{D}\del\zeta)}{\e+\zeta} 
-\e\,\frac{\alpha\del\zeta\cdot(\vec{D}\del\zeta)}{(\e+\zeta)^2} 
\right.
\\
\left.+\e \del\cdot(\vec{D}\del P) 
+ \frac{\alpha\f\cdot \del\zeta}{\e+\zeta}
+ \f\cdot\del P + \phi - L \right] d\x .
\end{multline}
Let us proceed term by term and let us assume that $\rho$ is bounded in $B_R$ 
uniformly as $\e\to 0$. This is a reasonable assumption, since we are assuming 
that the point $\x=\vec{0}$ is a repeller and hence does not belong to the 
attractor of the deterministic system on which we require the bound. Since $P$ 
and $\phi$ are continuous, 
\begin{multline}
\bigg\vert \int\limits_{B_R}  \rho \left[ \e \del\cdot(\vec{D}\del P) + 
\f\cdot\del P + \phi - L \right] d\x \bigg\vert 
\leq
\\
\frac{4\pi}{3}\,C^3\e^{3\eta} \, \max_{B_R}\big\lbrace \rho \abs{ \e 
\del\cdot(\vec{D}\del P) + \f\cdot\del P + \phi - L} \big\rbrace 
\end{multline}
so
\begin{equation}
\label{E:Estimate1}
%\boxed{
\lim_{\e\to 0}\int\limits_{B_R}  \rho \left[ \e \del\cdot(\vec{D}\del P) + 
\f\cdot\del P + \phi - L \right] d\x =0.
%}
\end{equation}

To study the other terms, we switch to polar coordinates, $(x_1,...,x_{n}) \to 
(r,\theta_1,...,\theta_{n-1})$ where $r\in[0,R]$, $\theta_1,...,\theta_{n-2} \in 
[0,\pi]$, $\theta_{n-1}\in[0,2\pi]$ and
\begin{equation}
d\x = r^{n-1} \sin^{n-2}(\theta_{1})...\sin(\theta_{n-2}) \, dr d\theta_1 ... 
d\theta_{n-1}.
\end{equation}
The quantity $\del\cdot(\vec{D}\del\zeta)$ is a fixed real number. Moreover, since $\zeta$ is a homogeneous, positive definite quadratic form of $\x$ and 
$\vec{D}$ is positive semi-definite (recall that 
$\vec{D}=\boldsymbol{\sigma}\boldsymbol{\sigma}^T$) we can write
\begin{gather}
\zeta(\x) = r^2 F(\theta_1,...,\theta_{n-1}) 	\\
\del\zeta\cdot(\vec{D}\del\zeta) = r^2 G(\theta_1,...,\theta_{n-1})
\end{gather}
for some strictly positive function $F$ and non-negative function $G$. 

Let
\begin{align}
F^* &= \min_{\theta_1,...,\theta_{n-1}} F(\theta_1,...,\theta_{n-1}),	\\
G^* &= \max_{\theta_1,...,\theta_{n-1}} G(\theta_1,...,\theta_{n-1}),
\end{align}
and
\begin{equation}
I = \int_{r=0}^{R} 
\int_{\theta_{n-1}=0}^{2\pi}\int_{\theta_{n-2}=0}^{\pi}...\int_{\theta_{1}=0}^{\pi} 
\frac{r^{n-1} \sin^{n-2}(\theta_{1})...\sin(\theta_{n-2})}{\e + r^2 
F(\theta_1,...,\theta_{n-1})} \, d\theta_1 ... d\theta_{n-1} dr.
\end{equation}
Then, we have
\begin{align}
\bigg\vert \int\limits_{B_R} \frac{\alpha \e \,\rho \, 
\del\cdot(\vec{D}\del\zeta)}{\e+\zeta} d\x  \,\bigg\vert 
&\leq \e \,\abs{\alpha \del\cdot(\vec{D}\del\zeta)} \,\max_{B_R}(\rho) \, I
\\[10pt]
&\leq \e \,\abs{\alpha \del\cdot(\vec{D}\del\zeta)}  \,\max_{B_R}(\rho)  
\,2\,\pi^{n-1}\int_{r=0}^{R} \frac{r^{n-1}}{\e+r^2F^*} dr.
\notag
\end{align}
If $n=2$, the last term can be integrated to give
\begin{equation}
\bigg\vert \int\limits_{B_R} \frac{\alpha \e \,\rho \, 
\del\cdot(\vec{D}\del\zeta)}{\e+\zeta} d\x  \,\bigg\vert   \leq \e \, 
\bigg\lbrace \abs{\alpha \del\cdot(\vec{D}\del\zeta)} \,\max_{B_R}(\rho)  
\,\frac{\pi^{n-1}}{F^*} \log\left(1+C^2\,F^*\e^{2\eta-1}\right) \bigg\rbrace
\end{equation}
while when $n\geq 3$ we can estimate
\begin{equation}
\bigg\vert \int\limits_{B_R} \frac{\alpha \e \,\rho \, 
\del\cdot(\vec{D}\del\zeta)}{\e+\zeta} d\x  \,\bigg\vert  
\leq 
\e^{1+\eta} \bigg\lbrace 2C\pi^{n-1}\abs{\alpha \del\cdot(\vec{D}\del\zeta)} 
\max_{B_R}(\rho) \max_{r\in[0,R]}\left(\frac{r^{n-1}}{\e+r^2F^*}\right) 
\bigg\rbrace.
\end{equation}
It can be verified that the maximum of the last term is achieved at the endpoint 
$r=R=C\e^{\eta}$. Taking the limit $\e\to 0$ shows that for all $n\geq 2$
\begin{equation}
\label{E:Estimate2}
%\boxed{
\lim_{\e\to 0}\int\limits_{B_R} \frac{\alpha \e \,\rho \, 
\del\cdot(\vec{D}\del\zeta)}{\e+\zeta} d\x  =0.
%}
\end{equation}
Similarly,  we can show
\begin{multline}
\bigg\vert \int\limits_{B_R}  \frac{\alpha \e \rho 
\del\zeta\cdot(\vec{D}\del\zeta)}{(\e+\zeta)^2} d\x \,\bigg\vert \leq \e \, 
\abs{\alpha} 2\pi^{n-1} \max_{B_R}(\rho) 
\int_{r=0}^{R} \frac{r^{n+1}G^*}{\e^2 + r^4F^{*2}} dr 
\\[10pt]
\leq \begin{cases}
\displaystyle  \e \left[ \frac{\abs{\alpha} \pi^{n-1} G^* 
}{2F^{*2}}\,\max_{B_R}(\rho) \,\log( 1+C^4\,F^{*2}\e^{4\eta-2})\right], & n=2 
\\[10pt]
\displaystyle  \e^{1+\eta}\,\left[ 2C \abs{\alpha}\pi^{n-1} 
G^*\,\max_{B_R}(\rho) \,\max_{r\in[0,R]}\left( \frac{r^{n+1}}{\e^2 + 
r^4F^{*2}}\right)\right] , & n\geq 3
\end{cases}
\end{multline}
where, again, the last maximum is achieved at $r=R=C{\e}^{\eta}$. We conclude that
\begin{equation}
\label{E:Estimate3}
%\boxed{
\lim_{\e\to 0}\int\limits_{B_R}  \frac{\alpha \e \rho 
\del\zeta\cdot(\vec{D}\del\zeta)}{(\e+\zeta)^2} d\x   =0.
%}
\end{equation}
Finally, since $\f(0)=0$ and $\del\zeta$ is linear, the term $\f\cdot\del\zeta$ 
is a polynomial of $\x$ that only contains monomials of degree 2 and higher. 
Consequently, we can write 
\begin{equation}
\f\cdot\del\zeta = \sum_{m=1}^{\deg(\f)} r^{1+m} H_m(\theta_1,...,\theta_{n-1})
\end{equation}
for some continuous functions $H_m$ whose usual $L^{\infty}$ norm 
$\norm{H_m}_{\infty}$ is finite. Each term in this series can be considered 
separately; for each $m$ we have
\begin{align}
\bigg\vert \int\limits_{B_R} \rho \frac{\alpha r^{1+m} H_m}{\e+\zeta} d\x 
\bigg\vert 
&\leq 2\pi^{n-1}\abs{\alpha}\,\norm{H_m}_{\infty}\,\max_{B_R}(\rho) \int_{0}^{R} 
\frac{r^{n+m}}{\e + r^2F^*} dr	\\
&\leq
\e^{\eta} \bigg\lbrace 2C 
\pi^{n-1}\abs{\alpha}\,\norm{H_m}_{\infty}\,\max_{B_R}(\rho) \,\max_{r\in[0,R]} 
\left( \frac{r^{n+m}}{\e + r^2F^*} \right)
\bigg\rbrace, \notag
\end{align}
which tends to 0 as $\e\to 0$ (since $n\geq 2$, $m\geq 1$, and the last maximum 
is obtained at $r=R=C{\e}^{\eta}$). We therefore conclude that
\begin{equation}
\label{E:Estimate4}
%\boxed{
\lim_{\e\to 0}\int\limits_{B_R} \rho \frac{\alpha \f\cdot \del\zeta}{\e+\zeta} 
d\x =0.
%}
\end{equation}
Combining~\eqref{E:Estimate1}, \eqref{E:Estimate2}, \eqref{E:Estimate3}  and \eqref{E:Estimate4} 
proves~\eqref{E:VanishingNoise_Thesis} and, 
consequently,~\eqref{E:VanishingNoise_Thesis0}. Note that the proof presented 
here is valid for systems of dimension $n\geq2$ or higher; one-dimensional 
systems are not of interest as bounds can be computed analytically.

\section{Construction of $\zeta$}
\label{S:ConstructZeta}
If the unstable fixed point is a repeller and $\vec{J}_0$ can be diagonalized, a 
quadratic form $\zeta$ that satisfies~\eqref{E:ReducedBalanceCondition} can be 
constructed from the eigenvectors of $\vec{J}_0$. Let $\vec{U}$ 
denote the matrix of eigenvectors of $\vec{J}_0$ and $\boldsymbol{\Lambda}$ be 
the usual diagonal matrix of eigenvalues, such that
\begin{equation}
\vec{U}^{-1}\vec{J}_0 \vec{U} = \boldsymbol{\Lambda}.
\end{equation}
Note that $\boldsymbol{\Lambda}$ is a positive definite matrix since the 
unstable point is repelling. Letting $\vec{w} = \vec{U}^{-1} \x$, an appropriate 
choice of $\zeta$ is
\begin{equation}
\label{E:ExampleZeta}
\zeta =  \vec{w}^T\vec{w} = \x^T [\vec{U}^{-1}]^T \vec{U}^{-1} \x.
\end{equation}
In fact, we have
\begin{equation}
\begin{aligned}
\tilde{\f}\cdot\del\zeta
&= 2 \x^T \vec{J}_0^T \, [\vec{U}^{-1}]^T \vec{U}^{-1} \x\\
&= 2\x^T  [ \vec{U}\vec{U}^{-1}]^T\, \vec{J}_0^T \, [\vec{U}^{-1}]^T 
\vec{U}^{-1} \x  \\
&= 2 \vec{w}^T [ \vec{U}^{-1} \vec{J}_0 \,\vec{U}]^T \vec{w} \\
&= 2 \vec{w}^T \boldsymbol{\Lambda}\vec{w},
\end{aligned}
\end{equation}
which is positive for $\vec{x}\neq\vec{0}$ because $\boldsymbol{\Lambda}$ is 
positive definite. We conclude that~\eqref{E:ReducedBalanceCondition} holds for 
$\alpha>0$. Note, however, that this is not the only possible choice of $\zeta$.

%%%%%%%%%%%%%%%%%%%%%%%% 
%    ACKNOWLEDGEMENTS
%%%%%%%%%%%%%%%%
\section{Acknowledgements}
The authors thank G.\ Chini, C.\ R.\ Doering, and A.\ Wynn for their helpful 
comments and discussions. D.G.\ and S.C.\ are grateful for the hospitality of 
the Institute for Pure and Applied Mathematics, where this work was begun, and 
G.F.\ and D.G.\ are grateful for the hospitality of the Geophysical Fluid 
Dynamics program at the Woods Hole Oceanographic Institution, where most of this 
work was carried out. During part or all of this work, G.F. was supported by 
Imperial College London under the IC PhD Scholarship Resfs G82059, D.G.\ was 
supported by the US NSF Mathematical Physics award PHY-1205219, and D.H.\ and 
S.C. received funding from EPSRC under the grant EP/J011126/1 and support in 
kind from Airbus Operation Ltd., ETH Zurich (Automatic Control Laboratory), 
University of Michigan (Department of Mathematics), and University of 
California, Santa Barbara (Department of Mechanical Engineering). 
%

%%%%%%%%%%%%%%%%%%%%%%%% 
%    REFERENCES
%%%%%%%%%%%%%%%%%%%%%%%%
\small
\setlength\bibsep{0pt}
\bibliographystyle{./hplain}

%%%%%%%%%%%%%%%%%%%%%%%%%%%%%%%%%%%%%%%%%%%%%%%%%%%%%%%%%%%%%%%%%%%%%%%%%%%%%%%%%%%%%%%%%%%%
\end{document}